\documentclass[a4paper,twoside]{article}

\usepackage{projective}

\author{Mireille Boutin and Gregor Kemper}

\title{On Reconstructing Configurations of Points in $\PP^2$ from a
Joint Distribution of Invariants}

\date{April 5, 2004}

\begin{document}

\maketitle

\begin{abstract}
  Consider the diagonal action of the projective group $\PGL_3$ on~$n$
  copies of $\PP^2$. In addition, consider the action of the symmetric
  group $\Sigma_n$ by permuting the copies. In this paper we find a
  set of generators for the invariant field of the combined group
  $\Sigma_n \times \PGL_3$. As the main application, we obtain a
  reconstruction principle for point configurations in $\PP^2$ from
  their sub-configurations of five points. Finally, we address the
  question of how such reconstruction principles pass down to
  subgroups.
\end{abstract}


\section*{Introduction} \label{0sIntro}

Consider the problem of recognizing a flat object from its shadow.
This is a common problem in computer vision where one often represents
objects by the boundary of their image on a picture.  For simplicity,
assume that the flat object is represented by a finite set of points
$p_1,\ldots,p_n\in \RR^3$.  Rotations and translations of such a flat
object in $\RR^3$ (almost always) induce a transformation of the image
points $P_1,\ldots,P_n\in \RR^2$ which can be written as
\begin{equation} \label{0eqAction}
  P_i \mapsto \frac{\left( \begin{array}{cc}
        a_{11} & a_{12} \\
        a_{21} & a_{22}
      \end{array}\right)
    P_i + \left(\begin{array}{c}
        a_{13} \\
        a_{23}
      \end{array}  \right)}{ (a_{31},a_{32}) P_i+a_{33} }, \text{ for
    all }i=1,\ldots,n,
  \text{ with }
  \left( \begin{array}{ccc}
      a_{11} & a_{12} & a_{13} \\
      a_{21} & a_{22} & a_{23} \\
      a_{31} & a_{32} & a_{33}
    \end{array}\right) \in GL(3)
\end{equation}
(where of course we have to assume that the above denominator does not
vanish). In the computer vision community, this group action is called
the projective group action ($PGL_3(\RR)=GL_3(\RR)/\RR^*$) and plays
an important role in many applications.

In order to be able to recognize a flat object from its shadow, we
thus need to be able to determine whether two sets of $n$ points in
the plane lie in the same orbit under the simultaneous action of the
projective group on each of the points.  More precisely, given
$P_1,\ldots,P_n\in\RR^2$ and $Q_1,\ldots,Q_n\in\RR^2$, we need to be
able to determine whether there exists a projective transformation
$g\in PGL_3(\RR)$ such that $g(P_i)=Q_i$, for all $i=1,\ldots,n$.
However, in many applications, the point correspondence between the
two objects is unknown: a priori, we ignore which point is going to be
mapped to which.  So, more generally, given any $P_1,\ldots,P_n$ and
$Q_1,\ldots,Q_n\in\RR^2$, we need to be able to determine whether
there exists a permutation $\pi\in \Sigma_n$ and a projective
transformation $g\in PGL_3(\RR)$ such that $g(P_i)=Q_{\pi(i)}$, for
all $i=1,\ldots,n$.

In an earlier publication~[\citenumber{Boutin.Kemper}], we considered
the analogue problem with the Euclidean group $\AO(2)$, which is a
subgroup of the projective group.  More precisely, we considered those
projective transformations whose matrix is given by
\[
\left( \begin{array}{ccc}
    a_{11} & a_{12} & a_{13} \\
    a_{21} & a_{22} & a_{23} \\
    a_{31} & a_{32} & a_{33}
  \end{array}\right) =
\left( \begin{array}{ccc}
    a_{11} & a_{12} & a_{13} \\
    a_{21} & a_{22} & a_{23} \\
    0 & 0 & 1
  \end{array}\right), \text{ with }\left( \begin{array}{cc}
    a_{11} & a_{12}  \\
    a_{21} & a_{22}
  \end{array}\right)\in \Or(2),
\]
where $\Or$ denotes the orthogonal group.

The squared distances $d_{i,j}=<P_i-P_j, P_i-P_j >$ are invariants
under the Euclidean group action, i.e.~they remain unchanged when
$P_i$ and $P_j$ are replaced by $g(P_i)$ and $g(P_j)$ respectively,
for any $g\in \AO(2)$.  Given $P_1,\ldots,P_n\in\RR^2$ and
$Q_1,\ldots,Q_n\in\RR^2$, it is a well known fact that $<P_i-P_j,
P_i-P_j >=<Q_i-Q_j, Q_i-Q_j >$ for every $i,j=1,\ldots,n$ if and only
if there exists a Euclidean transformation mapping $P_i$ to $Q_i$, for
every $i=1,\ldots, n$.  In order to take care of the labeling
ambiguity, we have tried to compare the {\em distribution} of the
pairwise distances of each point configurations, i.e.~the number of
times each value of the distances occurs.  Although there exist
$P_1,\ldots,P_n\in \RR^2$ and $Q_1,\ldots,Q_n\in \RR^2$ which have the
same distribution of distances but are not the same up to a relabeling
of the point and a Euclidean transformation, such examples are fairly
rare.  In fact, we have shown that there exists a non-zero polynomial
$f$ in $2n$ variables such that if $f(P_1,\ldots,P_n)\neq 0$, then the
point configuration $P_1,\ldots,P_n$ is uniquely determined up to a
Euclidean transformation and a relabeling.  In other words, if
$f(P_1,\ldots,P_n)\neq 0$, then for any $Q_1,\ldots,Q_n\in\RR^2$ with
the same distribution of distances as $P_1,\ldots,P_n$, there exists a
relabeling $\pi\in \Sigma_n$ and a Euclidean transformation $g\in
\AO(2)$ such that $g( P_i)=Q_{\pi(i)}$, for every $i=1,\ldots,n$.

In~[\citenumber{Boutin.Kemper}], we also considered the group of area
preserving affine transformations, which consists of those matrices
\[
\left( \begin{array}{ccc}
    a_{11} & a_{12} & a_{13} \\
    a_{21} & a_{22} & a_{23} \\
    a_{31} & a_{32} & a_{33}
  \end{array}\right) =
\left( \begin{array}{ccc}
    a_{11} & a_{12} & a_{13} \\
    a_{21} & a_{22} & a_{23} \\
    0 & 0 & 1
  \end{array}\right) \text{ with } \left| \begin{array}{cc}
    a_{11} & a_{12}  \\
    a_{21} & a_{22}
  \end{array}\right| =\pm 1.
\]
In that case, we looked at the distribution of the triangular areas
$\Delta_{i_1 i_2 i_3}= \frac{1}{2} \left| \left(
    P_{i_2}-P_{i_1}\right) \right. \times \left. \left(
    P_{i_3}-P_{i_1}\right) \right|$, for every distinct $i_1, i_2,
i_3\in \{1,\ldots,n \}$.  Obviously, areas remain unchanged under any
area-preserving affine transformation.  In a similar manner as with
the Euclidean group, we were able to show that there exists a non-zero
polynomial $f$ in $2n$ variables such that if $f(P_1,\ldots,P_n)\neq
0$, then $P_1,\ldots,P_n$ is uniquely determined, up to a relabeling
and an area preserving linear transformation, by the distribution of
its triangular areas.  In other words, there exists a Zariski-open set
of point configurations $(P_1,\ldots,P_n)\in \left(\RR^2\right)^n$
which are completely determined, up to an area-preserving affine
transformation and a relabeling, by the distribution of the triangular
areas between the $P_i's$.

We are now ready to attack the general case of a projective
transformation on $\PP^2(\RR)$.  In fact, everything we are about to
say holds for the more general case of the projective group $PGL_3(K)$
acting on the two-dimensional projective space $\PP^2(K)$ for any
infinite field $K$. In this context, the action given
in~\eqref{0eqAction} corresponds to the action on the subsets of
projective points of the form $(x:y:1)$.  In \sref{1sC}, we start by
obtaining a generating set of invariants for the diagonal (=
simultaneous) action of the projective group on~$n$ copies of
$\PP^2(K)$.  Some of these invariants turn out to be redundant and we
obtain a full set of relationships between them. These relationships
will be used over and over in the following. In classical invariant
theory, a theorem giving a full generating set of invariants of some
group $G$ acting diagonally on~$n$ copies of the natural
representation is often called the {\em first fundamental theorem} for
that group $G$, and then a theorem giving all relations between the
generators is called the {\em second fundamental theorem}. The
generating set we give has already appeared in \mycite{Olver}, but the
determination of the relations is, to the best of our knowledge, new.
In \sref{2sN5}, we consider the case $n = 5$ and take the action of
the symmetric group $\Sigma_5$ into account. We find two
invariants~$a$ and~$b$ which generate the invariant field
$K\left(\left(\PP^2(K)\right)^5\right)^{\Sigma_5 \times \PGL_3}$. This
is a crucial step toward the case of general $n$, which we attack in
\sref{3sGen}. In that section, we find a generating set of the field
of invariants of $\Sigma_n$ and $\PGL_3(K)$. In particular, given
$P_1,\ldots,P_n\in\PP^2(K)$, we consider the joint distribution of the
$a$'s and $b$'s evaluated at every
$P_{i_1},P_{i_2},P_{i_3},P_{i_4},P_{i_5}\in \{ P_1,\ldots,P_n\}$, with
$i_1,i_2,i_3,i_4,i_5$ distinct.  The final result
(\cref{3cSeparating}) of \sref{3sGen} states that there exists a
Zariski open subset $\Omega$ of $\left( \PP^2(K)\right)^n$ such that
any $(P_1,\ldots,P_n)\in \Omega$ is completely determined, up to a
projective transformation and a relabeling of the points, by the joint
distribution of the $a$'s and $b$'s.

Note that in this paper and in~[\citenumber{Boutin.Kemper}] the
general approach to reconstructing objects (modulo group actions) is
to consider the distribution of specified sub-objects (e.g.,
triangles, pentagons). In the final section of the paper, we formalize
and generalize this approach. Then we prove a theorem which under
rather mild hypotheses allows to transport this approach from one
group to an arbitrary subgroup. Combining this with
\cref{3cSeparating} and with the results
from~[\citenumber{Boutin.Kemper}], we obtain reconstruction theorems
for arbitrary subgroups of $\PGL_3$ and of area preserving
transformations.

\head{Acknowledgment} This research was initiated during a visit of
both authors at the Mathematical Sciences Research Institute in
Berkeley. We thank Michael Singer and Bernd Sturmfels for the
invitation.

\section{The first and second fundamental theorem for PGL$_3$} \label{1sC}

The main goal of this section is to prove what in classical invariant
theory would be termed the first and second fundamental theorem for
$\PGL_3$.

Let $K$ be any infinite field. We write $\PP^2 = \PP^2(K)$ for the
two-dimensional projective space, and $\PGL_3 = \PGL_3(K) =
\GL_3(K)/K^*$ for the projective group acting on $\PP^2$. Points from
$\PP^2$ are given by their homogeneous coordinates
$(\alpha_1:\alpha_2:\alpha_3)$ with $\alpha_i \in K$ not all zero. The
first lemma is an elementary fact from projective geometry.

\begin{lemma} \label{1lPGLOrbits}
  Let $P_1 \upto P_4 \in \PP^2$ be four projective points such that no
  three of them are collinear. Then there exists $g \in \PGL_3$ such
  that
  \[
  g(P_1) = (1:0:0), \quad g(P_2) = (0:1:0), \quad g(P_3) = (0:0:1),
  \quad \text{and} \quad g(P_4) = (1:1:1).
  \]
  This~$g$ is unique.
\end{lemma}

\begin{proof}
  For each point $P_i$ take a representative $v_i \in K^3$. Since
  $v_1$, $v_2$, and $v_3$ are linearly independent, we have
  \[
  v_3 = \alpha_1 v_1 + \alpha_2 v_2 + \alpha_3 v_3
  \]
  with $\alpha_i \in K$. Since no three of the $P_i$ are collinear,
  all $\alpha_i$ are non-zero. Thus we can choose the $v_i$ in such a
  way that $\alpha_i = 1$ for all~$i$. There exists a $\phi \in
  \GL_3(K)$ such that
  \[
  \phi(v_1) = (1,0,0), \quad \phi(v_2) = (0,1,0), \quad \text{and}
  \quad \phi(v_3) = (0,0,1).
  \]
  Now $v_4 = v_1 + v_2 + v_3$ implies $\phi(v_4) = (1,1,1)$. This
  proves the existence of $g \in \PGL_3$ with the claimed properties.

  To prove the uniqueness of~$g$ assume we have $\psi \in \GL_3(K)$
  with
  \[
  \psi(v_1) = \beta_1 \cdot (1,0,0), \quad \psi(v_2) = \beta_2 \cdot
  (0,1,0), \quad \psi(v_3) = \beta_3 \cdot (0,0,1), \quad \text{and}
  \quad \psi(v_4) = \beta_4 \cdot (1,1,1),
  \]
  where $\beta_i \in K \setminus \{0\}$ for all~$i$. Then $v_1 + v_2 +
  v_3 = v_4$ implies $(\beta_1,\beta_2,\beta_3) =
  (\beta_4,\beta_4,\beta_4)$, so all $\beta$'s are equal, and $\psi =
  \beta_1 \cdot \phi$. Therefore $\psi$ and $\phi$ define the same
  element in $\PGL_3$, which proves uniqueness.
\end{proof}

Following \mycite{Olver}, we describe rational invariants of~$n$
projective points. So let~$n$ be a positive integer and take $3 n$
indeterminates $x_{i,j}$ ($i \in \{1 \upto n\}$, $j \in \{0,1,2\}$).
We write $K(\underline{x})$ for the field of rational functions in the
$x_{i,j}$, and
\[
K(\underline{x})_0 := K\left(\frac{x_{1,1}}{x_{1,0}} \upto
  \frac{x_{n,1}}{x_{n,0}},\frac{x_{1,2}}{x_{1,0}} \upto
  \frac{x_{n,2}}{x_{n,0}}\right),
\]
which is the function field on $\left(\PP^2(K)\right)^n$.
Alternatively, $K(\underline{x})_0$ can be defined as the field of all
rational functions $f \in K(\underline{x})$ where for each $i \in \{1
\upto n\}$ the numerator and the denominator of~$f$ are homogeneous as
polynomials in $x_{i,0}$, $x_{i,1}$, $x_{i,2}$, and of the same
degree. We have a diagonal action of $\PGL_3$ on
$\left(\PP^2\right)^n$, which induces an action on the function field
$K(\underline{x})_0$ by $g(f) = f \circ g^{-1}$. For indices
$i_0,i_1,i_2 \in \{1 \upto n\}$ define the ``bracket''
\[
[i_0,i_1,i_2] := \det\left(x_{i_\nu,\mu}\right)_{\nu,\mu = 0,1,2} \in
K(\underline{x}),
\]
and for $i,j,k,l,m \in \{1 \upto n\}$ pairwise distinct define
\begin{equation} \label{1eqC}
  c_{i,j,k,l,m} := \frac{[i,j,k] [i,l,m]}{[i,j,l] [i,k,m]} \in
  K(\underline{x})_0.
\end{equation}
It is easy to see that the $c_{i,j,k,l,m}$ are $\PGL_3$-invariants. We
write
\[
K(\underline{x})_0^{\PGL_3(K)} = \{f \in K(\underline{x})_0 \mid g(f)
= f \ \text{for all} \ g \in \PGL_3(K)\}
\]
for the field of all $\PGL_3$-invariants. The first part of the
following theorem already appeared in \mycite{Olver} (though his
statement is slightly different).

\begin{theorem} \label{1tCs}
  With the above notation we have
  \begin{enumerate}
  \item (First fundamental theorem for $\PGL_3$.) The $c_{i,j,k,l,m}$
    generate the field of $\PGL_3$-invariants, i.e.,
    \[
    K(\underline{x})_0^{\PGL_3(K)} = K\left(c_{i,j,k,l,m} \mid
      i,j,k,l,m \in \{1 \upto n\} \ \text{pairwise distinct}\right).
    \]
  \item The $c_{i,j,k,l,m}$ separate $\PGL_3$-orbits on a dense open
    subset of $\left(\PP^2\right)^n$. More precisely, let $P_1 \upto
    P_n \in \PP^2(K)$ be points such that no three of them are
    collinear, and let $Q_1 \upto Q_n \in \PP^2(K)$ be further points
    such that
    \[
    c_{i,j,k,l,m}(P_1 \upto P_n) = c_{i,j,k,l,m}(Q_1 \upto Q_n)
    \]
    for all $i,j,k,l,m \in \{1 \upto n\}$ pairwise distinct (implying
    that no zero-division occurs when evaluating the $c_{i,j,k,l,m}$
    at $Q_1 \upto Q_n$), then there exists $g \in \PGL_3(K)$ such that
    \[
    g(P_i) = Q_i
    \]
    for all $i \in \{1 \upto n\}$.
  \end{enumerate}
\end{theorem}

\begin{proof}
  Let $d \in K[\underline{x}]$ be the product of all $[i,j,k]$ with $1
  \le i < j < k \le n$. For $P_1 \upto P_n \in \PP^2(K)$ with
  homogeneous coordinates $P_i = (\xi_{i,0}:\xi_{i,1},\xi_{i,2})$, we
  have that no three of the $P_i$ are collinear if and only if
  $d(\underline{\xi}) \ne 0$.

  We first treat the case $n \le 4$. By \lref{1lPGLOrbits}, all $(P_1
  \upto P_n) \in \left(\PP^2\right)^n$ where~$d$ takes a non-zero
  value lie in one single $\PGL_3$-orbit. Hence every invariant $f \in
  K(\underline{x})_0^{\PGL_3(K)}$ is constant on the set of all these
  $(P_1 \upto P_n)$. By \lref{1lDense} (which is proved after this
  lemma), $f$ is constant. This proves~(a) and~(b) of the lemma.

  Now assume $n \ge 5$ and consider the subset
  \[
  T := \left\{(P_1 \upto P_n) \in \left(\PP^2(K)\right)^n \mid P_1 =
    (1:0:0),P_2 = (0:1:0), P_3 = (0:0:1), P_4 = (1:1:1)\right\}
  \]
  of $\left(\PP^2(K)\right)^n$. \lref{1lPGLOrbits} implies that the
  set $\PGL_3\cdot T := \{g(\underline{P}) \mid g \in \PGL_3(K),
  (\underline{P}) \in T\}$ contains all $(P_1 \upto P_n) \in
  \left(\PP^2(K)\right)^n$ such that no three of $P_1$, $P_2$, $P_3$,
  and $P_4$ are collinear, so in particular $\PGL_3\cdot T$ contains
  all $(P_1 \upto P_n)$ where~$d$ takes a non-zero value. Thus
  \lref{1lDense} implies:
  \begin{equation} \label{1eqDense}
    \text{If two rational functions coincide on} \ \PGL_3\cdot T, \
    \text{they coincide as rational functions.}
  \end{equation}
  
  To prove~(a), take $0 \ne f \in K(\underline{x})_0^{\PGL_3(K)}$.
  Being a rational function in the $x_{i,j}$, $f$ can be written as $f
  = a/b$ with $a,b \in K\left[x_{i,j} \mid i \in \{1 \upto n\}, j \in
    \{0,1,2\}\right]$ coprime. It is easy to see that for each $i \in
  \{1 \upto n\}$, $a$ and~$b$ are homogeneous as polynomials in
  $x_{i,0}$, $x_{i,1}$, and $x_{i,2}$.
  Indeed, for $\alpha \in K$, let $\phi_{i,\alpha}$ be the
  $K$-automorphism of $K(\underline{x})$ which sends $x_{i,\nu}$ to
  $\alpha \cdot x_{i,\nu}$ and $x_{j,\nu}$ to itself for $j \ne
  i$. Then $f \in K(\underline{x})_0$ implies that
  $\phi_{i,\alpha}(a)/\phi_{i,\alpha}(b)= \phi_{i,\alpha}(f) = f =
  a/b$, so
  \[
  b \phi_{i,\alpha}(a) = \phi_{i,\alpha}(b) a.
  \]
  By the coprimality of~$a$ and~$b$ this implies that~$b$ divides
  $\phi_{i,\alpha}(b)$. Since $\phi_{i,\alpha}(b)$ and $b$ contain the
  same monomials, this means that $\phi_{i,\alpha}(b)$ is a scalar
  multiple of~$b$. Thus~$b$ is homogeneous as a polynomial in
  $x_{i,0}$, $x_{i,1}$, and $x_{i,2}$. The same argument works
  for~$a$.  Thus for a vector $(v_1 \upto v_n) \in \left(K^3 \setminus
  \{0\}\right)^n$, whether or not $b(v_1 \upto v_n)$ is 0 depends only
  on the class of $(v_1 \upto v_n)$ in $\left(\PP^2(K)\right)^n$.  So
  we can write $Z$ for the vanishing set of~$b$ {\em as a subset of}
  $\left(\PP^2(K)\right)^n$. By way of contradiction, assume that $T
  \subseteq Z$.  For $\phi \in \GL_3(K)$, the $\PGL_3$-invariance
  of~$f$ implies $\phi(a)/\phi(b) = a/b$, hence
  \[
  b \phi(a) = \phi(b) a.
  \]
  By the coprimality of~$a$ and~$b$ this implies that~$b$ divides
  $\phi(b)$, so if~$b$ vanishes at a point $(v_1 \upto v_n) \in
  \left(K^3\right)^n$, then~$b$ also vanishes at $\left(\phi^{-1}(v_1)
  \upto \phi^{-1}(v_n)\right)$ for all $\phi \in \GL_3(K)$.  Therefore
  the assumption $T \subseteq Z$ implies that $\PGL_3\cdot T \subseteq
  Z$. Now~\eqref{1eqDense} implies the contradiction $b = 0$.
  
  Having seen that $b$ does not vanish identically on $T$, we may
  define the restriction of~$f$ on $T$ and obtain a rational function
  on $T$:
  \[
  f|_T = F\left(\frac{x_{5,1}}{x_{5,0}} \upto
    \frac{x_{n,1}}{x_{n,0}},\frac{x_{5,2}}{x_{5,0}} \upto
    \frac{x_{n,2}}{x_{n,0}}\right),
  \]
  with $F$ a rational function in $2 (n - 4)$ arguments. Remembering
  the definition of the $c_{i,j,k,l,m}$ and evaluating them on $T$
  yields for $i > 4$:
  \begin{equation} \label{1eqCs}
    c_{3,2,4,i,1}|_T = x_{i,1}/x_{i,0} \quad \text{and} \quad
    c_{2,3,4,i,1}|_T = x_{i,2}/x_{i,0}.
  \end{equation}
  Hence
  \[
  f|_T = F\left(c_{3,2,4,5,1} \upto c_{3,2,4,n,1},c_{2,3,4,5,1} \upto
    c_{2,3,4,n,1}\right)|_T.
  \]
  So~$f$ and $F\left(c_{3,2,4,5,1} \upto c_{3,2,4,n,1},c_{2,3,4,5,1}
  \upto c_{2,3,4,n,1}\right)$ are two functions in
  $K(\underline{x})_0^{\PGL_3(K)}$ which coincide on $T$, hence they
  also coincide on $\PGL_3\cdot T$. Now~\eqref{1eqDense} implies that
  these functions coincide as elements of $K(\underline{x})$. This
  proves~(a).
  
  After these preparations, the proof of~(b) is easy. First, the
  hypothesis that none of the denominators vanish when evaluating the
  $c_{i,j,k,l,m}$ at $(Q_1 \upto Q_n)$ implies that, as for the $P_i$,
  no three of the $Q_i$ are collinear. Thus by \lref{1lPGLOrbits}
  there exist $\phi_1,\phi_2 \in \PGL_3(K)$ such that
  \[
  \left(\phi_1(P_1) \upto \phi_1(P_n)\right) \in T \quad \text{and}
  \quad \left(\phi_2(Q_1) \upto \phi_2(Q_n)\right) \in T.
  \]
  The hypothesis in~(b) and the invariance of the $c_{i,j,k,l,m}$
  imply that
  \[
  c_{i,j,k,l,m}\left(\phi_1(P_1) \upto \phi_1(P_n)\right) =
  c_{i,j,k,l,m}\left(\phi_2(Q_1) \upto \phi_2(Q_n)\right).
  \]
  Now~\eqref{1eqCs} implies that $\phi_1(P_i) = \phi_2(Q_i)$ for $i
  \ge 5$. But for $i \le 4$ this also holds by the definition of $T$.
  This completes the proof of~(b).
\end{proof}

The previous proof used the following elementary fact.

\begin{lemma} \label{1lDense}
  Let $f,g \in K(x_1 \upto x_m)$ be rational functions in~$m$
  indeterminates over the infinite field $K$, and let $h \in K[x_1
  \upto x_m] \setminus \{0\}$ be a non-zero polynomial. If
  \[
  f(\xi_1 \upto \xi_m) = g(\xi_1 \upto \xi_m)
  \]
  for all $\xi_1 \upto \xi_m \in K$ such that $h(\xi_1 \upto \xi_m)
  \ne 0$ and the evaluations of~$f$ and~$g$ at $(\xi_1 \upto \xi_m)$
  are defined, then $f = g$ (as rational functions).
\end{lemma}

\begin{proof}
  After subtracting~$g$ from~$f$ we may assume that $g = 0$. Next we
  multiply~$h$ by the denominator of~$f$, which does not change the
  hypothesis of the lemma. But now we can also multiply~$f$ by its
  denominator, so we may assume $f \in K[x_1 \upto x_m]$.
  We use induction on~$m$. By way of contradiction, assume that $f \ne
  0$.
  Since $K$ is infinite, there exists $\xi_m \in K$ such that $f_1 :=
  f(x_1 \upto x_{m-1},\xi_m) \in K[x_1 \upto x_{m-1}]$ is non-zero,
  and the same for $h_1 := h(x_1 \upto x_{m-1},\xi_m)$. If $m = 1$,
  this is an immediate contradiction to the hypothesis. If $m > 1$, we
  obtain a contradiction by induction.
\end{proof}

A major step in our argument is the study of relations between the
$c_{i,j,k,l,m}$. Let $P_0$ be a polynomial ring over $K$ with
indeterminates $C_{i,j,k,l,m}$ for $i,j,k,l,m \in \{1 \upto n\}$
pairwise distinct. Consider the homomorphism
\[
\mapl{\Phi}{P_0}{K(\underline{x})_0^{\PGL_3(K)}}%
{C_{i,j,k,l,m}}{c_{i,j,k,l,m}},
\]
and let $I_0 \subseteq P_0$ be the kernel of $\Phi$. Thus $I_0$ is the
ideal of relations between the $c_{i,j,k,l,m}$.

\begin{theorem}[Second fundamental theorem for $\PGL_3$] \label{1tSecond}
  The ideal $I_0$ is generated by the following relations:
  \begin{equation} \label{1eqEqual}
    \begin{aligned}
    C_{i,j,k,l,m} - C_{i,k,j,m,l},\\
    C_{i,j,k,l,m} - C_{i,l,m,j,k},\\
    C_{i,j,k,l,m} - C_{i,m,l,k,j},
    \end{aligned}
  \end{equation}
  \begin{equation} \label{1eqInverses}
    C_{i,j,k,l,m} \cdot C_{i,j,l,k,m} - 1,
  \end{equation}
  \begin{equation}  \label{1eqSum1}
    C_{i,j,k,l,m} + C_{i,j,m,l,k} - 1,
  \end{equation}
  \begin{equation}  \label{1eqTriad5}
    C_{i,j,k,l,m} - C_{m,j,k,l,i} \cdot C_{j,i,k,l,m},
  \end{equation}
  \begin{equation}  \label{1eqTriad6}
    C_{i,j,k,l,m} - C_{i,r,k,l,m} \cdot C_{i,j,k,l,r},
  \end{equation}
  where \eqref{1eqEqual}--\eqref{1eqTriad5} are for all $i,j,k,l,m \in
  \{1 \upto n\}$ pairwise distinct, and~\eqref{1eqTriad6} is for all
  $i,j,k,l,m,r \linebreak \in \{1 \upto n\}$ pairwise distinct.
\end{theorem}

\begin{proof}
  We first check that the relations given
  in~\eqref{1eqEqual}--\eqref{1eqTriad6} lie in $I_0$.
  For~\eqref{1eqEqual} and~\eqref{1eqInverses}, this is immediately
  seen from the definition of the
  $c_{i,j,k,l,m}$. For~\eqref{1eqSum1}, observe that
  \[
  c_{i,j,k,l,m} + c_{i,j,m,l,k} - 1 = \frac{[i,j,k] [i,l,m] - [i,j,m]
    [i,l,k] - [i,j,l] [i,k,m]}{[i,j,l] [i,k,m]}.
  \]
  The numerator is a function of five vectors $v_i,v_j,v_k,v_l,v_m \in
  K^3$. Fixing $v_i$, we see that the numerator is an alternating
  bilinear form in the arguments $v_j,v_k,v_l,v_m$. But an alternating
  bilinear form in four three-dimensional vectors has to be zero,
  hence the relation~\eqref{1eqSum1}. Next we check~\eqref{1eqTriad5}
  and~\eqref{1eqTriad6}:
  \[
  c_{m,j,k,l,i} \cdot c_{j,i,k,l,m} = \frac{[m,j,k] [m,l,i] [j,i,k]
    [j,l,m]}{[m,j,l] [m,k,i] [j,i,l] [j,k,m]} = c_{i,j,k,l,m},
  \]
  and
  \[
  c_{i,r,k,l,m} \cdot c_{i,j,k,l,r} = \frac{[i,r,k] [i,l,m] [i,j,k]
    [i,l,r]}{[i,r,l] [i,k,m] [i,j,l] [i,k,r]} = c_{i,j,k,l,m}.
  \]
  
  \newcommand{\bc}{\overline{C}}%
  
  Let $I \subseteq P_0$ be the ideal generated by the
  relations~\eqref{1eqEqual}--\eqref{1eqTriad6}, so $I \subseteq
  I_0$. We need to show the reverse inclusion $I_0 \subseteq I$. To
  this end, let $R := P_0/I$ be the residue class ring, and for
  $i,j,k,l,m \in \{1 \upto n\}$ distinct, write $\bc_{i,j,k,l,m} :=
  C_{i,j,k,l,m} + I \in R$ for the residue class of
  $C_{i,j,k,l,m}$. It follows from~\eqref{1eqInverses} that
  $\bc_{i,j,k,l,m}$ is invertible in $R$, i.e.
  \begin{equation} \label{1eqUnit}
    \bc_{i,j,k,l,m} \in R^\times,
  \end{equation}
  where $R^\times$ denotes the group of units in $R$. Consider the
  $K$-subalgebra $R_0 \subseteq R$ generated by all $\bc_{2,3,4,i,1}$
  and $\bc_{3,2,4,i,1}$ for $i \in \{5 \upto n\}$. Moreover, set $S :=
  R_0 \cap R^\times$ and
  \[
  R_1 := S^{-1} R_0 \subseteq R.
  \]
  We claim that $R_1 = R$, so we need to prove that for all $i,j,k,l,m
  \in \{1 \upto n\}$ pairwise distinct $\bc_{i,j,k,l,m}$ lies in
  $R_1$. For this purpose we first remark that if there exists a
  permutation $\pi$ of the set $\{j,k,l,m\}$ such that
  $\bc_{i,\pi(j),\pi(k),\pi(l),\pi(m)} \in R_1$, then also
  $\bc_{i,j,k,l,m} \in R_1$. Before giving the proof, we summarize the
  claim by stating
  \begin{equation} \label{1eqS4}
    \exists\ \pi \in \Sigma_{\{j,k,l,m\}}: \
    \bc_{i,\pi(j),\pi(k),\pi(l),\pi(m)} \in R_1 \quad \Rightarrow
    \quad \bc_{i,j,k,l,m} \in R_1,
  \end{equation}
  where $\Sigma$ denotes the symmetric group. Indeed, if~$\pi$ is the
  permutation given by $j \mapsto k$, $k \mapsto j$, $l \mapsto m$,
  and $m \mapsto l$, then~\eqref{1eqS4} follows directly
  from~\eqref{1eqEqual}. The same is true if~$\pi$ exchanges~$j$
  with~$l$ and~$k$ with~$m$. Furthermore, if~$\pi$ is given by $j
  \mapsto j$, $k \mapsto m$, $l \mapsto l$, and $m \mapsto k$, then
  \[
  \bc_{i,j,k,l,m} =  1 -  \bc_{i,j,m,l,k} = 1 -
  \bc_{i,\pi(j),\pi(k),\pi(l),\pi(m)}
  \]
  by~\eqref{1eqSum1}, so~\eqref{1eqS4} holds for this~$\pi$,
  too. Finally, if~$\pi$ is given by $j \mapsto j$, $k \mapsto l$, $l
  \mapsto k$, and $m \mapsto m$, then
  \[
  \bc_{i,j,k,l,m} =  \bc_{i,j,l,k,m}^{-1} =
  \bc_{i,\pi(j),\pi(k),\pi(l),\pi(m)}^{-1}
  \]
  by~\eqref{1eqInverses}. If $\bc_{i,\pi(j),\pi(k),\pi(l),\pi(m)} \in
  R_1$, then $\bc_{i,\pi(j),\pi(k),\pi(l),\pi(m)} = f/g$ with $f \in
  R_0$ and $g \in S$, so
  \[
  f = g \cdot \bc_{i,\pi(j),\pi(k),\pi(l),\pi(m)} \in R^\times \cap
  R_0 = S
  \]
  by~\eqref{1eqUnit}. Thus $\bc_{i,j,k,l,m} = g/f \in R_1$,
  so~\eqref{1eqS4} holds for this~$\pi$, too. But the four particular
  $\pi$'s considered so far generate the symmetric group
  $\Sigma_{\{j,k,l,m\}}$, so~\eqref{1eqS4} follows in general.
  
  Now we prove that $\bc_{i,j,k,l,m} \in R_1$ for all $i,j,k,l,m \in
  \{1 \upto n\}$ distinct. If $2 \notin \{i,j,k,l,m\}$, then
  \[
  \bc_{i,j,k,l,m} = \bc_{i,2,k,l,m} \cdot \bc_{i,j,k,l,2}
  \]
  by~\eqref{1eqTriad6}. Thus we are done if we can show that all
  $\bc_{i,j,k,l,m}$ with $2 \in \{i,j,k,l,m\}$ lie in $R_1$. In other
  words, we may assume that $2 \in \{i,j,k,l,m\}$. By~\eqref{1eqS4} we
  may even assume that $i = 2$ or $k = 2$. Furthermore, if $3 \notin
  \{i,j,k,l,m\}$, then
  \[
  \bc_{i,j,k,l,m} = \bc_{i,3,k,l,m} \cdot \bc_{i,j,k,l,3}
  \]
  by~\eqref{1eqTriad6}, so we may assume that $3 \in \{i,j,k,l,m\}$
  (preserving $i = 2$ or $k = 2$). If $i \ne 2$ and $i \ne 3$, we may
  assume $j = 2$ and $m = 3$ by using~\eqref{1eqS4}. Then
  \[
  \bc_{i,j,k,l,m} = \bc_{i,2,k,l,3} = \bc_{3,2,k,l,i} \cdot
  \bc_{2,i,k,l,3}
  \]
  by~\eqref{1eqTriad5}. This means that we may assume $i = 2$ or $i =
  3$ and, moreover (using~\eqref{1eqS4}) that $\{i,k\} = \{2,3\}$. Now
  if $1 \notin \{j,k,l,m\}$, then
  \[
  \bc_{i,j,k,l,m} = \bc_{i,1,k,l,m} \cdot \bc_{i,j,k,l,1}
  \]
  by~\eqref{1eqTriad5}, so after all we may assume $1 \in
  \{j,k,l,m\}$. Using~\eqref{1eqS4} again, we can achieve $\{i,k,l\} =
  \{1,2,3\}$. If $4 \notin \{j,m\}$ then
  \[
  \bc_{i,j,k,l,m} = \bc_{i,4,k,l,m} \cdot \bc_{i,j,k,l,4}
  \]
  by~\eqref{1eqTriad6}, so we may assume $4 \in
  \{j,k,l,m\}$. Using~\eqref{1eqS4} again, we obtain $(i,j,k,l,m) =
  (2,3,4,l,1)$ or $(i,j,k,l,m) = (3,2,4,l,1)$. Thus $\bc_{i,j,k,l,m}$
  lies in $R_1$ as claimed, which completes the proof that $R_1 = R$.
  
  We still need to prove that $I_0 \subseteq I$, so take $f \in
  I_0$. We have $f + I \in R = R_1$, so there exist polynomials $g,h
  \in P$ involving only the indeterminates $C_{2,3,4,i,1}$ and
  $C_{3,2,4,i,1}$ ($i \in \{5 \upto n\}$) such that $h f - g \in I$
  and $h + I \in R^\times$. Since $f \in I_0$ we obtain
  \begin{equation} \label{1eqG0}
    \Phi(g) = \Phi(h) \cdot \Phi(f) = 0.
  \end{equation}
  Consider the subset $T_0 \in \left(\PP^2(K)\right)^n$ consisting of
  all those $(P_1 \upto P_n) \in \left(\PP^2(K)\right)^n$ with $P_1 =
  (1:0:0)$, $P_2 = (0:1:0)$, $P_3 = (0:0:1)$, and $P_4 = (1:1:1)$,
  such that the first coordinate of all $P_i$ with $i \ge 5$ is
  non-zero (see the proof of \tref{1tCs}). Let $\bar{g} \in
  K\left(x_{5,2}/x_{5,0} \upto x_{n,2}/x_{n,0},x_{5,1}/x_{5,0} \right.
  \upto \linebreak \left. x_{n,1}/x_{n,0}\right)$ be the rational
  function obtained from~$g$ by substituting each $C_{2,3,4,i,1}$ by
  $x_{i,2}/x_{i,0}$ and each $C_{3,2,4,i,1}$ by $x_{i,1}/x_{i,0}$ (for
  $i \ge 5$). It follows from~\eqref{1eqCs} that~$\bar{g}$ and
  $\Phi(g)$ coincide as functions on $T_0$. Thus by~\eqref{1eqG0},
  $\bar{g}$ vanishes on $T_0$. But a rational function in
  indeterminates $x_{i,2}/x_{i,0}$ and $x_{i,1}/x_{i,0}$ ($i \ge 5$)
  vanishing on $T_0$ must be zero by~\lref{1lDense}. Since the
  $x_{i,2}/x_{i,0}$ and $x_{i,1}/x_{i,0}$ are algebraically
  independent it follows that $g = 0$. Now $h f - g \in I$ implies $h
  f \in I$. Together with $h + I \in R^\times$, this implies $f \in
  I$. This completes the proof that $I_0 \subseteq I$.
\end{proof}

\begin{remark*}
  One can see from the proof of \tref{1tSecond} that the invariant
  field $K(\underline{x})_0^{\PGL_3(K)}$ is in fact generated by the
  $c_{2,3,4,i,1}$ and $c_{3,2,4,i,1}$ with $i \in \{5 \upto n\}$
  (other choices are possible), and that these $2 (n - 4)$ generators
  are algebraically independent. So in particular
  $K(\underline{x})_0^{\PGL_3(K)}$ is purely transcendental over
  $K$. The ``extended'' generating system containing all
  $c_{i,j,k,l,m}$ is nevertheless more suitable for our purposes,
  since it is permuted by the action of the symmetric group $\Sigma_n$
  on the indices of each $c_{i,j,k,l,m}$.
\end{remark*}

In Sections~\ref{2sN5} and~\ref{3sGen} we will need the following
lemma, which gives some ``non-relations''.

\begin{lemma} \label{1lRelations}
  We keep the notation of \tref{1tSecond}.
  \begin{enumerate}
  \item The relations given in~\eqref{1eqEqual} are the only
    equalities that exist between the $c_{i,j,k,l,m}$. More precisely,
    if
    \[
    c_{i',j',k',l',m'} = c_{i,j,k,l,m},
    \]
    then $i' = i$, and the list $[j',k',l',m']$ is one of $[j,k,l,m]$,
    $[k,j,m,l]$, $[l,m,j,k]$, or $[m,l,k,j]$.
  \item For each $\nu \in \{1,2,3\}$, let
    $i_\nu,j_\nu,k_\nu,l_\nu,m_\nu \in \{1 \upto n\}$ be pairwise
    distinct indices, and suppose that
    \begin{equation} \label{1eqTriadc}
      c_{i_1,j_1,k_1,l_1,m_1} = c_{i_2,j_2,k_2,l_2,m_2} \cdot
      c_{i_3,j_3,k_3,l_3,m_3}.
    \end{equation}
    Then
    \[
    s :=
    |\{i_1,j_1,k_1,l_1,m_1,i_2,j_2,k_2,l_2,m_2,i_3,j_3,k_3,l_3,m_3\}|
    \in \{5,6\},
    \]
    i.e., only five or six indices occur in the above relation.
  \item If $s = 6$ in~(b), then $i_1 = i_2 = i_3$.
  \item If $i_1 = i_2 = i_3$ does {\em not} hold in~(b), then
    \[
    \{i_2,i_3\} = \{j_1,m_1\} \quad \text{or} \quad \{i_2,i_3\} =
    \{k_1,l_1\}.
    \]
  \end{enumerate}
\end{lemma}

\begin{proof}
  To prove~(a), assume $c_{i',j',k',l',m'} = c_{i,j,k,l,m}$. Then
  every bracket $[\nu,\mu,\eta]$ occurring in $c_{i',j',k',l',m'}$
  must contain the index~$i$, hence $i' = i$. Moreover,
  $c_{i',j',k',l',m'}$ must have the bracket $[i,j,k]$ or $[i,k,j]$ in
  its numerator and bracket $[i,j,l]$ or $[i,l,j]$ in its denominator,
  hence the claim.
  
  Now assume the hypothesis of~(b). First observe that if some
  index~$\nu$ occurs in this relation, it must occur at least twice,
  since otherwise one side of~\eqref{1eqTriadc} would involve the
  indeterminates $x_{\nu,\mu}$ while the other side would not.
  
  We will study the behavior of both sides of~\eqref{1eqTriadc} when
  we equate some of the arguments $v_i$. More precisely, for $i,j \in
  \{1 \upto n\}$ distinct and for $f \in K[\underline{x}]$ an
  irreducible polynomial, set $w_{\{i,j\}} := 1$ if~$f$ lies in the
  ideal generated by $x_{i,0} - x_{j,0}$, $x_{i,1} - x_{j,1}$, and
  $x_{i,2} - x_{j,2}$, and set $w_{\{i,j\}} := 0$ otherwise. Extend
  $w_{\{i,j\}}$ to a function $\left(K(\underline{x}) \setminus
    \{0\}\right) \to \ZZ$ by using the rule $w_{\{i,j\}}(f g) =
  w_{\{i,j\}}(f) + w_{\{i,j\}}(g)$. Thus for $i,j,k,l,m \in \{1 \upto
  n\}$ pairwise distinct we have:
  \begin{align*}
    & w_{\{j,k\}}(c_{i,j,k,l,m}) = w_{\{l,m\}}(c_{i,j,k,l,m}) = 1, \\
    & w_{\{j,l\}}(c_{i,j,k,l,m}) = w_{\{k,m\}}(c_{i,j,k,l,m}) = -1,
    \quad \text{and} \\
    & w_{\{\nu,\mu\}}(c_{i,j,k,l,m}) = 0 \quad \text{for} \quad
    \{\nu,\mu\} \notin
    \left\{\{j,k\},\{l,m\},\{j,l\},\{k,m\}\strut\right\}.
  \end{align*}
  These equations will be used frequently in the sequel.
  \eref{1eqTriadc} implies
  \[
  w_{\{j_1,k_1\}}(c_{i_2,j_2,k_2,l_2,m_2}) +
  w_{\{j_1,k_1\}}(c_{i_3,j_3,k_3,l_3,m_3}) = 1,
  \]
  so $\{j_1,k_1\} \in
  \left\{\{j_2,k_2\},\{l_2,m_2\},\{j_3,k_3\},\{l_3,m_3\}\strut\right\}$.
  Possibly exchanging factors on the left hand side
  of~\eqref{1eqTriadc} (which does not change any of the assertions of
  part~(b), (c), or (d) of the lemma), we may assume that $\{j_1,k_1\}
  = \{j_2,k_2\}$ or $\{j_1,k_1\} = \{l_2,m_2\}$.
  Using~\eqref{1eqEqual}, we may now reorder the indices
  $j_2,k_2,l_2,m_2$ in such a way that
  \[
  j_1 = j_2 \quad \text{and (consequently)} \quad k_1 = k_2.
  \]
  Using the same argument with $w_{\{l_1,m_1\}}$ yields $\{l_1,m_1\}
  \in \left\{\{l_2,m_2\},\{j_3,k_3\},\{l_3,m_3\}\strut\right\}$.
  
  First consider the case $\{l_1,m_1\} = \{l_2,m_2\}$.
  Then~\eqref{1eqTriadc} becomes
  \[
  c_{i_1,j_1,k_1,l_1,m_1} = c_{i_2,j_1,k_1,l_1,m_1} \cdot
  c_{i_3,j_3,k_3,l_3,m_3} \quad \text{or} \quad
  c_{i_1,j_1,k_1,l_1,m_1} = c_{i_2,j_1,k_1,m_1,l_1} \cdot
  c_{i_3,j_3,k_3,l_3,m_3}.
  \]
  It follows that $\{j_3,k_3,l_3,m_3\} = \{j_1,k_1,l_1,m_1\}$, since
  otherwise some $w_{\{\nu,\mu\}}$ would take the value~1 on the right
  hand side of the above equation but~0 on the left hand side. (For
  example, if $j_3 \notin \{j_1,k_1,l_1,m_1\}$, this would apply to
  $w_{\{j_3,k_3\}}$.) But then $i_1 = i_2 = i_3$, since otherwise some
  $i_\nu$ would occur only once as an index in~\eqref{1eqTriadc},
  which cannot happen. Thus $s = 5$ and we are done with proving
  (b)--(d) in this case. (In fact, carrying the arguments further
  shows that this case cannot occur.)
  
  It remains to consider the cases $\{l_1,m_1\} = \{j_3,k_3\}$ or
  $\{l_1,m_1\} = \{l_3,m_3\}$. As above we may use~\eqref{1eqEqual} to
  reorder the indices $j_3,k_3,l_3,m_3$ in such a way that $l_1 = l_3$
  and (consequently) $m_1 = m_3$, so~\eqref{1eqTriadc} becomes
  \[
  c_{i_1,j_1,k_1,l_1,m_1} = c_{i_2,j_1,k_1,l_2,m_2} \cdot
  c_{i_3,j_3,k_3,l_1,m_1}.
  \]
  Since $w_{\{j_1,l_1\}}(c_{i_1,j_1,k_1,l_1,m_1}) = -1$, we must have
  $\{j_1,l_1\} \in
  \left\{\{j_1,l_2\},\{k_1,m_2\},\{j_3,l_1\},\{k_3,m_1\}\strut\right\}$.
  The second and the fourth possibilities would violate the
  distinctness of $i_2,j_1,k_1,l_2,m_2$ or $i_3,j_3,k_3, \linebreak
  l_1,m_1$, respectively, so we have $l_1 = l_2$ or $j_1 = j_3$.

  Consider the case $l_1 = l_2$. Then
  \[
  c_{i_1,j_1,k_1,l_1,m_1} = c_{i_2,j_1,k_1,l_1,m_2} \cdot
  c_{i_3,j_3,k_3,l_1,m_1}.
  \]
  Applying the above argument again (using $w_{\{k_1,m_1\}}$) shows
  $m_1 = m_2$ or $k_1 = k_3$. But if $m_1 = m_2$, then
  $w_{\{l_1,m_1\}}$ takes different values on the different sides of
  the above equation, so $k_1 = k_3$. Thus
  \[
  c_{i_1,j_1,k_1,l_1,m_1} = c_{i_2,j_1,k_1,l_1,m_2} \cdot
  c_{i_3,j_3,k_1,l_1,m_1}.
  \]
  We have $m_1 \ne m_2$ and $j_1 = j_2 \ne m_2$, hence
  \[
  0 = w_{\{l_1,m_2\}}(c_{i_1,j_1,k_1,l_1,m_1}) = 1 +
  w_{\{l_1,m_2\}}(c_{i_3,j_3,k_1,l_1,m_1}),
  \]
  implying $j_3 = m_2$. We are left with
  \begin{equation} \label{1eqCase1}
    c_{i_1,j_1,k_1,l_1,m_1} = c_{i_2,j_1,k_1,l_1,m_2} \cdot
    c_{i_3,m_2,k_1,l_1,m_1}.
  \end{equation}
  The set $T := \{j_1,k_1,l_1,m_1,m_2\}$ has~5 (distinct) elements.
  Assume $i_1 \notin T$. Then for the prime polynomial $[i_1,j_1,k_1]$
  to appear in the numerator of $c_{i_2,j_1,k_1,l_1,m_2} \cdot
  c_{i_3,m_2,k_1,l_1,m_1}$ we must have $i_2 = i_1$. Furthermore, $i_3
  = i_1$, since otherwise $[i_1,l_1,m_1]$ could not appear in that
  numerator. Now assume $i_2 \notin T$. Then $i_1 = i_2$, since
  otherwise $[i_2,j_1,k_1]$ would appear only once on the right hand
  side of~\eqref{1eqCase1} and not at all on the left hand side.
  Moreover, $[i_2,k_1,m_2]$ does not appear on the left hand side, so
  it must be cancelled on the right hand side, so $i_3 = i_2$.
  Likewise, if $i_3 \notin T$, then $[i_3,m_2,k_1]$ must be cancelled
  on the right hand side of~\eqref{1eqCase1}, so $i_2 = i_3$, and
  $[i_3,k_1,m_1]$ must appear on the left hand side, so $i_1 = i_3$.
  Thus we have seen that if any of the $i_\nu$ lie in $T$, then $i_1 =
  i_2 = i_3$ and thus $s = 6$. The other possibility is that all
  $i_\nu$ lie in $T$. But then $s = 5$ and $i_1 = m_2$ (otherwise the
  indices on the left hand side of~\eqref{1eqCase1} would not be
  distinct), $i_2 = m_1$, and $i_3 = j_1$. So we are in one of the
  cases described by part~(d) of the lemma. Thus parts~(b)--(d) are
  proved in the case $l_1 = l_2$.
  
  Now consider the remaining case $j_1 = j_3$. We have
  \[
  c_{i_1,j_1,k_1,l_1,m_1} = c_{i_2,j_1,k_1,l_2,m_2} \cdot
  c_{i_3,j_1,k_3,l_1,m_1}.
  \]
  Considering $w_{\{k_1,m_1\}}$ yields $m_1 = m_2$ or $k_1 = k_3$. The
  possibility $k_1 = k_3$ is ruled out by considering
  $w_{\{j_1,k_1\}}$, so $m_1 = m_2$. Thus $l_1 \ne l_2$ (since $l_1 =
  l_2$ was considered above) and $k_1 = k_2 \ne l_2$ yield
  \[
  0 = w_{\{l_2,m_1\}}(c_{i_1,j_1,k_1,l_1,m_1}) = 1 +
  w_{\{l_2,m_1\}}(c_{i_3,j_1,k_3,l_1,m_1}),
  \]
  so $k_3 = l_2$, and we obtain
  \[
  c_{i_1,j_1,k_1,l_1,m_1} = c_{i_2,j_1,k_1,l_2,m_1} \cdot
  c_{i_3,j_1,l_2,l_1,m_1}.
  \]
  In this case we consider the set $T := \{j_1,k_1,l_1,m_1,l_2\}$ of
  size~5. Using exactly the same arguments as in the previous case, we
  conclude that either $s = 6$ and $i_1 = i_2 = i_3$, or $s = 5$ and
  $i_1 = l_2$, $i_2 = l_1$, and $i_3 = k_1$. So parts~(b)--(d) of the
  lemma are proved in this case, too.
\end{proof}

Before we go on, it is useful to introduce some notation which
deviates slightly from the notation introduced before
\tref{1tSecond}. By~\eqref{1eqEqual} and \lref{1lRelations}(a) there
are precisely $n (n-1) (n-2) (n-3) (n-4)/4$ distinct
$c_{i,j,k,l,m}$. We take as many indeterminates as follows: For
$i,j,k,l,m \in \{1 \upto n\}$ pairwise distinct with $j =
\min\{j,k,l,m\}$ let $C_{i,j,k,l,m}$ be an indeterminate over $K$. For
$i,j,k,l,m$ distinct but {\em not} meeting the additional constraint
that $j = \min\{j,k,l,m\}$, we define $C_{i,j,k,l,m}$ by imposing the
equations
\begin{equation} \label{1eqEqualC}
  C_{i,j,k,l,m} = C_{i,k,j,m,l} = C_{i,l,m,j,k} = C_{i,m,l,k,j},
\end{equation}
which reflect~\eqref{1eqEqual}. Let $P$ be the polynomial ring
generated by the $C_{i,j,k,l,m}$, and let $I$ be the kernel of the
homomorphism $P \to K(\underline{x})$ of $K$-algebras sending
$C_{i,j,k,l,m}$ to $c_{i,j,k,l,m}$. Thus $I$ is the ideal of relations
between the $c$'s. The distinction between the polynomial rings $P_0$
(introduced before \tref{1tSecond}) and $P$ may seem a bit subtle, but
introducing $P$ ultimately renders our notation much simpler. $P_0$
will not be used anymore in the sequel.


\section{The case $n = 5$} \label{2sN5}

In this section we will work out a set of generating invariants for
$K\left(\underline{x}\right)_0^{\Sigma_n \times \PGL_3(K)}$ in the
case $n = 5$. Here and in the sequel we write $\Sigma_n$ for the
symmetric group in~$n$ symbols. Recall that $K$ is an infinite field
and $K\left(\underline{x}\right)_0 = K\left(x_{i,j}/x_{i,0} \mid i = 1
  \upto n, j = 1,2\right)$ is the function field of
$\left(\PP^2(K)\right)^n$. We will also use the $\PGL_3(K)$-invariants
$c_{i,j,k,l,m}$ defined in~\eqref{1eqC} and the indeterminates
$C_{i,j,k,l,m}$ defined at the end of \sref{1sC}. \tref{1tSecond} and
\lref{1lRelations} give information on the ideal $I$ of relations
between the $c_{i,j,k,l,m}$. As we consider the case $n = 5$, there
are precisely 30 $C_{i,j,k,l,m}$. We denote the group of all
permutations of these 30 elements by $\Sigma_{30}$. (Note that this is
a slight deviation from the notation $\Sigma_n$ for the symmetric
group in~$n$ symbols.) Any such permutation acts on the polynomial
ring $P$ generated by the $C_{i,j,k,l,m}$. The crucial step in this
section is the proof of the following lemma.

\begin{lemma} \label{2lPerm}
  Let $\phi \in \Sigma_{30}$ be a permutation of the $C_{i,j,k,l,m}$
  with $\{i,j,k,l,m\} = \{1,2,3,4,5\}$ and $j =
  \min\{j,k,l,m\}$. Assume that $\phi$ maps the ideal $I \subset P$
  into itself. Then there exists a permutation $\pi \in \Sigma_5$ of
  the numbers $1 \upto 5$ such that for all indices $i,j,k,l,m$ we
  have
  \[
  \phi(C_{i,j,k,l,m}) = C_{\pi(i),\pi(j),\pi(k),\pi(l),\pi(m)}.
  \]
\end{lemma}

\begin{proof}
  Take $i,j,k,l,m \in \{1,2,3,4,5\}$ pairwise distinct with $j =
  \min\{j,k,l,m\}$ (meaning $j = 1$ if $i \ne 1$ and $j = 2$
  otherwise), and suppose
  \[
  \phi(C_{i,j,k,l,m}) = C_{r,s,t,u,v}
  \]
  with $\{r,s,t,u,v\} = \{1,2,3,4,5\}$, $s =
  \min\{s,t,u,v\}$. By~\eqref{1eqSum1} we have $C_{i,j,k,l,m} +
  C_{i,j,m,l,k} - 1 \in I$, hence our hypothesis implies
  $C_{r,s,t,u,v} + \phi(C_{i,j,m,l,k}) - 1 \in I$. On the other hand,
  $C_{r,s,t,u,v} + C_{r,s,v,u,t} - 1 \in I$, so $\phi(C_{i,j,m,l,k}) -
  C_{r,s,v,u,t} \in I$. By \lref{1lRelations}(a) it follows that
  \begin{equation} \label{2eqPerm1}
    \phi(C_{i,j,m,l,k}) = C_{r,s,v,u,t}.
  \end{equation}
  Using~\eqref{1eqInverses} we see that $C_{i,j,k,l,m} C_{i,j,l,k,m} -
  1 \in I$, hence $C_{r,s,t,u,v} \phi(C_{i,j,l,k,m}) - 1 \in I$. But
  also $C_{r,s,t,u,v} C_{r,s,u,t,v} - 1 \in I$, and the uniqueness of
  inverses in any ring (here: $P/I$) leads to
  \begin{equation} \label{2eqPerm2}
    \phi(C_{i,j,l,k,m}) = C_{r,s,u,t,v}.
  \end{equation}
  Repeated application of~\eqref{2eqPerm1} and~\eqref{2eqPerm2} shows
  that $\phi(C_{i,j,*,*,*}) = C_{r,s,+,+,+}$, where the $*$'s are the
  indices $k,l,m$ appearing in some order, and the $+$'s stand for
  $t,u,v$ appearing in the corresponding order.

  Now define a map $\map{\pi}{\{1 \upto 5\}}{\{1 \upto 5\}}$ as
  follows: For $i \in \{1 \upto 5\}$ there are unique $j,k,l,m$ with
  $\{i,j,k,l,m\} = \{1 \upto 5\}$ and $j < k < l < m$. Let $\pi(i)$ be
  the first index of $\phi(C_{i,j,k,l,m})$. There are precisely~6
  (distinct) $C_{r,s,t,u,v}$ with $r = \pi(i)$. By the above
  observation it follows that all these $C_{r,s,t,u,v}$ are images of
  suitable $C_{i,j,*,*,*}$ under~$\phi$. Therefore the hypothesis that
  $\phi$ permutes the $C$'s implies that $\pi$ is actually a
  permutation of the set $\{1 \upto 5\}$. Define $\phi_\pi$ by
  \[
  \phi_\pi(C_{i,j,k,l,m}) = C_{\pi(i),\pi(j),\pi(k),\pi(l),\pi(m)}.
  \]
  We wish to show that $\phi = \phi_\pi$. We know that for any
  $i,j,k,l,m \in \{1 \upto 5\}$ distinct we have $\phi(C_{i,j,k,l,m})
  = C_{\pi(i),s,t,u,v}$ with $s,t,u,v \in \{1 \upto 5\} \setminus
  \{\pi(i)\}$. Hence $(\phi_\pi^{-1} \circ \phi)(C_{i,j,k,l,m}) =
  C_{i,\pi^{-1}(s),\pi^{-1}(t),\pi^{-1}(u),\pi^{-1}(v)}$.
  Using~\eqref{1eqEqualC} we may reorder
  $\pi^{-1}(s),\pi^{-1}(t),\pi^{-1}(u),\phi^{-1}(v)$ in such a way
  that~$j$ (which has to be among the
  $\pi^{-1}(s),\pi^{-1}(t),\pi^{-1}(u),\phi^{-1}(v)$ since
  $\{i,j,k,l,m\} = \{1 \upto 5\} =
  \{i,\pi^{-1}(s),\pi^{-1}(t),\pi^{-1}(u),\pi^{-1}(v)\}$) appears
  first. Thus for $\{i,j,k,l,m\} = \{1 \upto 5\}$ we have
  \begin{equation} \label{2eqPi}
    \left(\phi_\pi^{-1} \circ \phi\right)(C_{i,j,k,l,m}) =
    C_{i,j,r,s,t}
  \end{equation}
  with $\{r,s,t\} = \{k,l,m\}$. The permutation $\phi_\pi$ sends the
  relation ideal $I \subset P$ to itself. In fact, $\pi$ induces an
  automorphism $\psi_\pi$ of $K(\underline{x})$ given by
  $\psi_\pi(x_{i,\nu}) = x_{\pi(i),\nu}$. With
  $\map{\Phi}{P}{K(\underline{x})}$ given by $\Phi(C_{i,j,k,l,m}) =
  c_{i,j,k,l,m}$, we clearly have $\Phi\left(\phi_\pi(f)\right) =
  \psi_\pi\left(\Phi(f)\right)$ for $f \in P$, hence $f \in I$ implies
  $\Phi\left(\phi_\pi(f)\right) = \psi_\pi(0) = 0$, so indeed
  $\phi_\pi(f) \in I$. To simplify notation, we may thus replace
  $\phi$ by $\phi_\pi^{-1} \circ \phi$. Then~\eqref{2eqPi} leads to
  \begin{equation} \label{2eqPhi}
    \phi(C_{i,j,k,l,m}) = C_{i,j,r,s,t} \quad \text{with} \quad
    \{r,s,t\} = \{k,l,m\},
  \end{equation}
  and we have to show that $\phi = \id$. By~\eqref{1eqTriad5} we have
  $C_{i,j,k,l,m} - C_{m,j,k,l,i} C_{j,i,k,l,m} \in I$. Since $\phi(I)
  \subseteq I$, this implies $C_{i,j,r,s,t} - C_{m,j,*,*,*}
  C_{j,i,*,*,*} \in I$ with the $*$'s standing for appropriate (as yet
  unknown) indices. By \lref{1lRelations}(d) we conclude $\{j,m\} =
  \{j,t\}$ or $\{j,m\} = \{r,s\}$. The second possibility is ruled out
  by the distinctness of $i,j,r,s,t$, hence $t = m$.

  Furthermore, we have $C_{i,j,k,l,m} + C_{i,j,m,l,k} - 1 \in I$
  by~\eqref{1eqSum1}, hence $C_{i,j,r,s,t} + \phi(C_{i,j,m,l,k}) - 1
  \in I$. On the other hand, $C_{i,j,r,s,t} + C_{i,j,t,s,r} - 1 \in
  I$, implying $\phi(C_{i,j,m,l,k}) =
  C_{i,j,t,s,r}$. Using~\eqref{1eqTriad5} again, we obtain
  $C_{i,j,m,l,k} - C_{k,j,m,l,i} C_{j,i,m,l,k} \in I$, hence
  $C_{i,j,t,s,r} - C_{k,j,*,*,*} C_{j,i,*,*,*} \in
  I$. \lref{1lRelations}(d) tells us that $\{k,j\} = \{j,r\}$ or
  $\{k,j\} = \{s,t\}$, hence $r = k$. Having seen that $t = m$ and $r
  = k$, we conclude that also $s = l$, since $\{r,s,t\} = \{k,l,m\}$
  by~\eqref{2eqPhi}. Thus~\eqref{2eqPhi} becomes $\phi(C_{i,j,k,l,m})
  = C_{i,j,k,l,m}$, so indeed $\phi = \id$. This completes the proof.
\end{proof}

Now we are ready to prove the main result of this section, which gives
a generating set for $K(\underline{x})_0^{\Sigma_n \times \PGL_3(K)}$
in the case $n = 5$.

\begin{theorem} \label{2tN5}
  With the $\PGL_3(K)$-invariants $c_{i,j,k,l,m}$ defined as
  in~\eqref{1eqC}, form
  \[
  a := \sum_{\pi \in \Sigma_5}
  c_{\pi(1),\pi(2),\pi(3),\pi(4),\pi(5)}^2 \quad \text{and} \quad b :=
  \sum_{\pi \in \Sigma_5} c_{\pi(1),\pi(2),\pi(3),\pi(4),\pi(5)}^4.
  \]
  If the characteristic of $K$ is not~2, then
  \[
  K(\underline{x})_0^{\Sigma_5 \times \PGL_3(K)} = K(a,b).
  \]
\end{theorem}

\begin{proof}
  Clearly~$a$ and~$b$ are invariant under $\Sigma_5$ and $\PGL_3(K)$,
  which shows the inclusion ``$\supseteq$''. So we need to prove
  ``$\subseteq$''. \tref{1tCs}(a) and~\eqref{1eqEqual} tell us that
  \[
  N := K(\underline{x})_0^{\PGL_3(K)} = K\left(c_{i,j,k,l,m} \mid
  \{i,j,k,l,m\} = \{1 \upto 5\},j = \min\{j,k,l,m\}\right).
  \]
  With $T$ an additional indeterminate, form the polynomial
  \[
  F := \prod
  \begin{Sb}
    \{i,j,k,l,m\} = \{1 \upto 5\}, \\ j = \min\{j,k,l,m\}
  \end{Sb}
  \left(T - c_{i,j,k,l,m}\right) \in N[T].
  \]
  A fairly easy computation using the computer algebra system
  Magma~[\citenumber{magma}] shows that the coefficients of $F$ lie in
  $L := K(a,b)$. In fact, using the relations given in
  \eqref{1eqEqual}--\eqref{1eqTriad5}, one can express all
  $c_{i,j,k,l,m}$ with $\{i,j,k,l,m\} = \{1 \upto 5\}$ and $j =
  \min\{j,k,l,m\}$ as rational functions in $c_{1,2,3,4,5}$ and
  $c_{2,1,3,4,5}$. So we get $c_{i,j,k,l,m} =
  f_{i,j,k,l,m}(c_{1,2,3,4,5},c_{2,1,3,4,5})$ with $f_{i,j,k,l,m} \in
  K(X,Y)$, where $X$ and $Y$ are new indeterminates. Details on what
  the $f_{i,j,k,l,m}$ actually are can be found in \rref{2rN5}(a). Now
  all that we need to do is express the 30 elementary symmetric
  functions of the $f_{i,j,k,l,m}$ in terms of the sum of squares and
  the sum of fourth powers of the $f_{i,j,k,l,m}$. Our Magma
  computation, which only involves rational functions in $X$ and $Y$,
  shows that (thanks to the special form of the $f_{i,j,k,l,m}$) this
  is indeed possible. It is in this computation that $\ch(K) \ne 2$ is
  required.

  We conclude that $N$ is the splitting field of $F$ over
  $L$. \lref{1lRelations}(a) implies that $F$ is a separable
  polynomial, hence $N$ is Galois as a field extension of $L$. By
  Galois theory we are done if we can show that the Galois group of
  $N$ over $L$ is contained in $\Sigma_5$ (in which case it will be
  equal to $\Sigma_5$). So take $\phi \in \Gal(N/L)$. Since $N$ is the
  splitting field of $F$, $\phi$ permutes the roots $c_{i,j,k,l,m}$ of
  $F$, and $\phi$ is determined by its permutation action on these
  roots. But since $\phi$ is a field automorphism, it preserves all
  relations between the $c_{i,j,k,l,m}$. This means that $\phi$,
  viewed as a permutation of the indeterminates $C_{i,j,k,l,m}$, maps
  the relation ideal $I$ into itself. Now it follows from
  \lref{2lPerm} that indeed $\phi \in \Sigma_5$. This completes the
  proof.
\end{proof}

\begin{remark} \label{2rN5}
  \begin{enumerate}
  \item It looks as if we had to evaluate all 30 of the
    $c_{i,j,k,l,m}$ with $\{i,j,k,l,m\} = \{1 \upto 5\}$ and $j =
    \min\{j,k,l,m\}$ in order to obtain the values of the
    invariants~$a$ and~$b$ appearing in \tref{2tN5}. But in fact they
    can all be expressed in terms of $c_{1,2,3,4,5}$ and
    $c_{2,1,3,4,5}$. Let us explain how. Form the set ${\mathcal M} :=
    \{X,Y,X/Y,(X-1)/(Y-1),X (1-Y)/(X-Y)\}$ with $X$ and $Y$
    indeterminates. For each $f \in {\mathcal M}$, also add $1/f$,
    $1-f$, $1/(1-f)$, $(f-1)/f$, and $f/(f-1)$ into $\mathcal M$, so
    that $\mathcal M$ contains a total of 30 rational functions in $X$
    and $Y$. Then all $c_{i,j,k,l,m}$ are obtained by substituting $X
    = c_{1,2,3,4,5}$ and $Y = c_{2,1,3,4,5}$ in the rational
    functions~$f$ from $\mathcal M$. This can be seen from the
    relations~\eqref{1eqEqual}--\eqref{1eqTriad5}. In particular,
    if we form
    \begin{equation} \label{2eqAB}
      A := \sum_{f \in \mathcal M} f^2 \quad \text{and} \quad B :=
      \sum_{f \in \mathcal M} f^4,
    \end{equation}
    we obtain
    \[
    a = 4 \cdot A(c_{1,2,3,4,5},c_{2,1,3,4,5}) \quad \text{and} \quad
    b = 4 \cdot B(c_{1,2,3,4,5},c_{2,1,3,4,5}).
    \]
  \item If follows from~\eqref{1eqCs} that $c_{3,2,4,5,1}$ and
    $c_{2,3,4,5,1}$ are algebraically independent over $K$. Thus the
    transcendence degree of $N := K(\underline{x})_0^{\PGL_3(K)} =
    K\left(c_{i,j,k,l,m} \mid \{i,j,k,l,m\} = \{1 \upto 5\}\right)$
    over $K$ is at least~2. But $L := K(\underline{x})_0^{\Sigma_5
    \times \PGL_3(K)} = N^{\Sigma_5}$ has the same transcendence
    degree, since $N$ is algebraic over $L$. From \tref{2tN5}, the
    transcendence degree of $L$ is at most~2. It follows that the
    transcendence degree of both fields is precisely~2, and the
    generating invariants $a_{1,2,3,4,5}$ and $b_{1,2,3,4,5}$ are
    algebraically independent. In particular, two is the smallest
    number of generating invariants for $L$ that we could have
    expected.
  \item We can also deal with the case $\ch(K) = 2$. In fact, this
    just requires a slight change of the invariants~$a$
    and~$b$. Instead of taking the sum of the squares and of fourth
    powers of the $c_{i,j,k,l,m}$, we need to take the second and
    fourth elementary symmetric functions in the $c_{i,j,k,l,m}$ with
    $\{i,j,k,l,m\} = \{1 \upto 5\}$ and $j = \min\{j,k,l,m\}$. In the
    context of part~(a) of this remark, we need to replace $A$ and $B$
    by the second and fourth elementary symmetric function in the
    $f$'s from $\mathcal M$. A Magma computation as mentioned in the
    proof of \tref{2tN5} then shows that {\em all} elementary
    symmetric functions in the $f$'s from $\mathcal M$ can be
    expressed as rational functions in $A$ and $B$.
  \item \mycite[page~141]{MLR98} determined a set of five
    $\PGL_3$-invariants of five points $(P_1 \upto P_5)$ which are
    also invariant under the action of the symmetric group $\Sigma_4$
    acting by permutations of the last four points
    $P_2,P_3,P_4,P_5$. These five invariants are permuted by the
    action of the complete permutation group $\Sigma_5$. The authors
    propose to take the values of these five invariants, ordered in
    increasing sequence, as invariants of $\Sigma_5 \times \PGL_3$.
  \end{enumerate}
\end{remark}

\section{The case of general~$n$} \label{3sGen}

In this section we attack the problem of finding generating invariants
of $K(\underline{x})_0^{\Sigma_n \times \PGL_3(K)}$ for a general
positive integer~$n$. Recall our notation. $K$ is an infinite field,
$K(\underline{x}) = K\left(x_{i,j} \mid 1 \le i \le n, \right.
\linebreak \left. 0 \le j \le 2\right)$ is a rational function field
in $3 n $ indeterminates over a field $K$, and for $i,j,k,l,m \in \{1
\upto n\}$ pairwise distinct we have a rational function
$c_{i,j,k,l,m}$ as given in~\eqref{1eqC}.  $P$ is a polynomial ring
over $K$ in $n (n-1) (n-2) (n-3) (n-4)/4$ indeterminates
$C_{i,j,k,l,m}$ labeled by $i,j,k,l,m \in \{1 \upto n\}$ pairwise
distinct with $j = \min\{j,k,l,m\}$.  Using~\eqref{1eqEqualC}, which
mirrors the equalities~\eqref{1eqEqual} existing between the
$c_{i,j,k,l,m}$, we define $C_{i,j,k,l,m}$ for {\em any} pairwise
distinct $i,j,k,l,m \in \{1 \upto n\}$. The ideal $I \subset P$ is the
kernel of the map $P \to K(\underline{x})$ sending each
$C_{i,j,k,l,m}$ to $c_{i,j,k,l,m}$; thus $I$ is the ideal of relations
of the $c_{i,j,k,l,m}$. The following lemma is tailored for proving
the main result, \tref{3tNgen}, of the section.

\begin{lemma} \label{3lPerm}
  Let $\psi$ be a permutation of the set
  \[
  {\mathcal M} := \left\{S \subseteq \{1 \upto n\} \mid |S| =
  5\right\},
  \]
  and for each $S \in \mathcal M$ let $\map{\pi_S}{S}{\psi(S)}$ be a
  bijection. Define a homomorphism $\map{\phi}{P}{P}$ of $K$-algebras by
  \[
  \phi(C_{i,j,k,l,m}) :=
  C_{\pi_S(i),\pi_S(j),\pi_S(k),\pi_S(l),\pi_S(m)} \quad \text{where}
  \quad S = \{i,j,k,l,m\}.
  \]
  (Note that $\phi$ is well-defined since the
  equalities~\eqref{1eqEqualC} are preserved.) If $\phi(I) \subseteq
  I$, then there exists a permutation $\pi \in \Sigma_n$ such that
  \[
  \phi(C_{i,j,k,l,m}) = C_{\pi(i),\pi(j),\pi(k),\pi(l),\pi(m)}
  \]
  for all~$i$, $j$, $k$, $l$, $m$.
\end{lemma}

\begin{proof}
  There is nothing to show for $n \le 5$, so we may assume $n \ge 6$.

  Let $T = \{i,j,k,l,m,r\} \subseteq \{1 \upto n\}$ be a set if six
  (distinct) elements. By~\eqref{1eqTriad6} we have $C_{i,j,k,l,m} -
  C_{i,r,k,l,m} \cdot C_{i,j,k,l,r} \in I$, hence by hypothesis also
  \begin{equation} \label{3eqTriad}
    \phi(C_{i,j,k,l,m}) - \phi(C_{i,r,k,l,m}) \cdot
  \phi(C_{i,j,k,l,r}) \in I.
  \end{equation}
  With $S := \{i,j,k,l,m\} \in \mathcal M$ we have
  $\phi(C_{i,j,k,l,m}) =
  C_{\pi_S(i),\pi_S(j),\pi_S(k),\pi_S(l),\pi_S(m)}$, and
  correspondingly for the other $C$'s occurring
  in~\eqref{3eqTriad}. Thus the union of all indices occurring
  in~\eqref{3eqTriad} is
  \[
  \tilde{T} := \psi\left(\{i,j,k,l,m\}\right) \cup
  \psi\left(\{i,r,k,l,m\}\right) \cup \psi\left(\{i,j,k,l,r\}\right).
  \]
  By \lref{1lRelations}(b), $\tilde{T}$ has at most six elements. On
  the other hand, the injectivity of $\psi$ implies that even the
  union of the $\psi$-images of just two different sets in $\mathcal
  M$ has at least six elements. Therefore
  \begin{equation} \label{3eqT}
    \tilde{T} = \psi\left(\{i,j,k,l,m\}\right) \cup
    \psi\left(\{i,r,k,l,m\}\right) = \psi\left(\{i,j,k,l,m\}\right)
    \cup \psi\left(\{i,j,k,l,r\}\right)
  \end{equation}
  and $\tilde{T}$ has precisely six elements. It follows that there
  exists $r' \in \tilde{T}$ such that $\psi\left(\{i,j,k,l,m\}\right)
  = \tilde{T} \setminus \{r'\}$. Likewise,
  $\psi\left(\{i,r,k,l,m\}\right) = \tilde{T} \setminus \{j'\}$ and
  $\psi\left(\{i,j,k,l,r\}\right) = \tilde{T} \setminus \{m'\}$ with
  $j',m' \in \tilde{T}$. Now we use the same argument with the roles
  of~$i$ and~$j$ interchanged. This yields
  \[
  \psi\left(\{j,i,k,l,m\}\right) \cup \psi\left(\{j,k,l,m,r\}\right) =
  \psi\left(\{j,i,k,l,m\}\right) \cup \psi\left(\{j,i,k,l,r\}\right)
  \]
  The second expression for $\tilde{T}$ in~\eqref{3eqT} is equal to
  the right hand side of above equation. Hence $\tilde{T} =
  \psi\left(\{j,i,k,l,m\}\right) \cup \psi\left(\{j,k,l,m,r\}\right)
  \supseteq \psi\left(\{j,k,l,m,r\}\right)$, so there exists $i' \in
  \tilde{T}$ with $\psi\left(\{j,k,l,m,r\}\right) = \tilde{T}
  \setminus \{i'\}$. In the same way, interchanging~$j$ and~$k$ yields
  $\psi\left(\{i,k,j,l,m\}\right) \cup \psi\left(\{i,j,l,m,r\}\right)
  = \tilde{T}$, so $\psi\left(\{i,j,l,m,r\}\right) = \tilde{T}
  \setminus \{k'\}$. Finally, interchanging~$j$ and~$l$ yields
  $\psi\left(\{i,l,k,j,m\}\right) \cup \psi\left(\{i,k,j,m,r\}\right)
  = \tilde{T}$, so $\psi\left(\{i,j,k,m,r\}\right) = \tilde{T}
  \setminus \{l'\}$. In summary, there exists a function
  $\map{\eta_T}{T}{\tilde{T} \subseteq \{1 \upto n\}}$ (which maps~$i$
  to~$i'$ etc.) such that $\psi\left(T \setminus \{\nu\}\right) =
  \tilde{T} \setminus \{\eta_T(\nu)\}$ for all $\nu \in T$. By
  hypothesis, $\psi$ is injective, so the same holds for $\eta_T$,
  hence $\eta_T(T) = \tilde{T}$. It follows that for any $S \in
  \mathcal M$ with $S \subset T$ we have
  \begin{equation} \label{3eqEta}
    \psi(S) = \eta_T(S),
  \end{equation}
  where the right hand side indicates element-wise application of
  $\eta_T$.

  It follows from~\eqref{3eqEta} that if two sets $S,S' \in \mathcal
  M$ have four elements in common, then also $\psi(S)$ and $\psi(S')$
  have four elements in common. In fact, $T := S \cup S'$ has six
  elements, hence $\psi(S) = \eta_T(S)$ and $\psi(S') =
  \eta_T(S')$. These are two subsets of size~5 inside the set
  $\eta_T(T)$ which has six elements, hence indeed $\psi(S)$ and
  $\psi(S')$ share four elements. Now take two subsets $T$,~$T'
  \subseteq \{1 \upto n\}$ with $|T| = |T'| = 6$ such that $S := T
  \cap T'$ has~5 elements. We will show that $\eta_T$ and $\eta_{T'}$
  coincide on $S$. Write
  \[
  T = S \cup \{j\} \quad \text{and} \quad T' = S \cup \{k\}
  \]
  with $j,k \in \{1 \upto n\}$. For $l \in S$ set $S_l := T' \setminus
  \{l\}$, so $S_l \in \mathcal M$. Then $|S_l \cap \left(T \setminus
  \{l\}\right)| = 4$ and $|S_l \cap S| = 4$, so, as noted above,
  $\psi(S_l)$ shares~4 elements with $\psi\left(T \setminus
  \{l\}\right) = \eta_T(T) \setminus \{\eta_T(l)\}$ and with $\psi(S)
  = \eta_T(S) = \eta_T(T) \setminus \{\eta_T(j)\}$. But $\psi(S_l)$
  cannot be a subset of $\eta_T(T)$ since this would imply
  \[
  \psi(S_l) = \eta_T\left(\eta_T^{-1}\left(\psi(S_l)\right)\right) =
  \psi\left(\eta_T^{-1}\left(\psi(S_l)\right)\right),
  \]
  contradicting the injectiveness of~$\psi$, since $S_l \not\subseteq
  T$. It follows that $\psi(S_l) = \eta_T\left(T \setminus
  \{j,l\}\right) \cup \{r_l\}$ with $r_l \in \{1 \upto n\} \setminus
  \eta_T(T)$. We can write this slightly simpler as $\psi(S_l) =
  \eta_T\left(S \setminus \{l\}\right) \cup \{r_l\}$. On the other
  hand, we have $S_l \subset T'$, so
  \[
  \psi(S_l) = \eta_{T'}(S_l) = \eta_{T'}\left(S \setminus \{l\}\right)
  \cup \{\eta_{T'}(k)\}.
  \]
  Intersecting the resulting equality $\eta_T\left(S \setminus
  \{l\}\right) \cup \{r_l\} = \eta_{T'}\left(S \setminus \{l\}\right)
  \cup \{\eta_{T'}(k)\}$ over all $l \in S$ yields $\bigcap_{l \in S}
  \{r_l\} = \{\eta_{T'}(k)\}$. Thus $r_l = \eta_{T'}(k)$ independently
  of~$l$, and $\eta_T\left(S \setminus \{l\}\right) = \eta_{T'}\left(S
  \setminus \{l\}\right)$ for all $l \in S$. This shows that
  $\eta_T(l) = \eta_{T'}(l)$ for all $l \in S$, as claimed.

  We proceed by taking any two subsets $T$,~$T' \subseteq \{1 \upto
  n\}$ with $|T| = |T'| = 6$. We can move from $T$ to $T'$ by
  successively exchanging elements. Using the above result, we see
  that $\eta_T$ and $\eta_{T'}$ coincide on $T \cap T'$. Thus we can
  define $\map{\pi}{\{1 \upto n\}}{\{1 \upto n\}}$ such that for every
  subset $T \subseteq \{1 \upto n\}$ with $|T| = 6$ the restriction
  $\pi|_{_T}$ coincides with $\eta_T$. Thus~\eqref{3eqEta} yields
  \[
  \psi(S) = \pi(S)
  \]
  for all $S \in \mathcal M$, where again the right hand side
  indicates element-wise application of~$\pi$. In particular, $\pi$ is
  injective, since otherwise $|\pi(S)| < 5$ for some $S \in \mathcal
  M$. Hence $\pi \in \Sigma_n$.

  Define $\map{\phi_\pi}{P}{P}$ by $\phi_\pi(C_{i,j,k,l,m}) :=
  C_{\pi(i),\pi(j),\pi(k),\pi(l),\pi(m)}$. We claim that $\phi =
  \phi_\pi$, which is equivalent to $\phi_{\pi^{-1}} \circ \phi =
  \id$. It is clear from the definition of $I$ that $\phi_\pi$ maps
  $I$ onto itself, hence $\left(\phi_{\pi^{-1}} \circ \phi\right)(I)
  \subseteq I$. For $i,j,k,l,m \in \{1 \upto n\}$ distinct we have
  \[
  \left(\phi_{\pi^{-1}} \circ \phi\right)(C_{i,j,k,l,m}) =
  C_{\pi^{-1}(\pi_S(i)),\pi^{-1}(\pi_S(j)),\pi^{-1}(\pi_S(k)),
  \pi^{-1}(\pi_S(l)),\pi^{-1}(\pi_S(m))},
  \]
  where $S := \{i,j,k,l,m\} \in \mathcal M$ and $\pi_S$ is given by
  the hypothesis of the lemma. Observe that $\left(\pi^{-1} \circ
  \pi_S\right)(S) = \pi^{-1}\left(\psi(S)\right) =
  \pi^{-1}\left(\pi(S)\right) = S$, so $\pi^{-1} \circ \pi_S$ is a
  bijection $S \to S$. Thus, in order to complete the proof, we may
  substitute $\phi$ by $\phi_{\pi^{-1}} \circ \phi$, and then we have
  the hypothesis that every $\pi_S$ is a bijection $S \to S$. Our goal
  is to show that all $\pi_S$ are equal to the identity.

  Assume that there exists an $S \in \mathcal M$ and an $i \in S$ such
  that $\pi(i) \ne i$. Set $j := \pi_S(i) \in S$ and write $S =
  \{i,j,k,l,m\}$. Moreover, choose any $r \in \{1 \upto n\} \setminus
  S$. By~\eqref{1eqTriad6} we have $C_{i,j,k,l,m} - C_{i,r,k,l,m}
  \cdot C_{i,j,k,l,r} \in I$. With $S' := \{i,r,k,l,m\}$ and $S'' :=
  \{i,j,k,l,r\}$ it follows that
  \[
  C_{\pi_S(i),\pi_S(j),\pi_S(k),\pi_S(l),\pi_S(m)} -
  C_{\pi_{S'}(i),\pi_{S'}(r),\pi_{S'}(k),\pi_{S'}(l),\pi_{S'}(m)}
  \cdot
  C_{\pi_{S''}(i),\pi_{S''}(j),\pi_{S''}(k),\pi_{S''}(l),\pi_{S''}(r)}
  \in I.
  \]
  By \lref{1lRelations}(c), this implies $\pi_S(i) = \pi_{S'}(i)$, but
  $\pi_S(i) = j \notin S' = \pi_{S'}(S')$. This contradiction shows
  that indeed all $\pi_S$ are the identity, completing the proof.
\end{proof}

To prove the main result of this section, we still need an elementary
lemma from field theory.

\begin{lemma} \label{3lGalois}
  Let $N = K(a_1 \upto a_m,b_1 \upto b_m)$ be a field extension of $K$
  generated by pairwise distinct elements $a_1 \upto a_m,b_1 \upto
  b_m$. Let $G \subseteq \Aut_K(N)$ be the group of all those
  $K$-automorphisms $\sigma$ of $N$ for which there exists $\pi \in
  \Sigma_m$ with $\sigma(a_i) = a_{\pi(i)}$ and $\sigma(b_i) =
  b_{\pi(i)}$ for all~$i$. Take indeterminates $X$, $T_1$, $T_2$, and
  consider the polynomial
  \[
  F := \prod_{i=1}^m \left(X - T_1 a_i - T_2 b_i\right) \in N[X,T_1,T_2].
  \]
  Let $L \subseteq N$ be the subextension generated by all
  coefficients of $F$. Then
  \[
  N^G = L.
  \]
\end{lemma}

\begin{proof}
  By the definition of $G$, any $\sigma \in G$ permutes the factors of
  $F$, hence $L \subseteq N^G$. We use Galois theory to prove the
  reverse inclusion. It follows from the construction of $L$ that
  $\prod_{i=1}^m (X - a_i)$ and $\prod_{i=1}^m (X - b_i)$ lie in
  $L[X]$. Therefore $N$ is the spitting field over $L$ of the
  polynomial $\prod_{i=1}^m \left((X - a_i) (X - b_i)\right)$. Hence
  $N/L$ is Galois, so $L = N^{\Gal(N/L)}$. If we can prove that
  $\Gal(N/L) \subseteq G$, then $N^G \subseteq N^{\Gal(N/L)} = L$, and
  we are done. So take any $\sigma \in \Gal(N/L)$. Writing $\sigma(F)$
  for the coefficient-wise application of $\sigma$ to $F$, we obtain
  \[
  \prod_{i=1}^m \left(X - T_1 a_i - T_2 b_i\right) = F = \sigma(F) =
  \prod_{i=1}^m \left(X - T_1 \sigma(a_i) - T_2 \sigma(b_i)\right).
  \]
  Since the zeros of a polynomial are uniquely determined up to
  permutations, there exists $\pi \in \Sigma_m$ such that $T_1
  \sigma(a_i) + T_2 \sigma(b_i) = T_1 a_{\pi(i)} + T_2 b_{\pi(i)}$ for
  all~$i$. It follows that indeed $\sigma \in G$.
\end{proof}

We can now give a generating set for the invariant field
$K(\underline{x})_0^{\Sigma_n \times \PGL_3(K)}$. We may assume $n \ge
5$, since for $n \le 4$ all invariants are constant (this is contained
in \tref{1tCs}(a)). Let $S \subseteq \{1 \upto n\}$ be a subset of
five elements. Set
\[
a_S:=\sum
\begin{Sb}
  i,j,k,l,m\ \text{with} \\
  \{1,j,k,l,m\} = S
\end{Sb}
c_{i,j,k,l,m}^2 \quad \text{and} \quad b_S:=\sum
\begin{Sb}
  i,j,k,l,m\ \text{with} \\
  \{1,j,k,l,m\} = S
\end{Sb}
c_{i,j,k,l,m}^4
\]
with the $c_{i,j,k,l,m}$ defined in~\eqref{1eqC}. These are clearly
functions in $K(\underline{x})$ which are invariant under the action
of $\PGL_3(K)$ and under all those permutations from $\Sigma_n$ which
map $S$ to itself.

\begin{theorem} \label{3tNgen}
  With the above notation, take additional indeterminates $X$, $T_1$,
  and $T_2$, assume the characteristic of $K$ is not~2, and form the
  polynomial
  \[
  F:=\prod
  \begin{Sb}
    S \subseteq \{1 \upto n\}, \\
    |S| = 5
  \end{Sb}
  \left(X - T_1 a_S - T_2 b_S\right) \in K(\underline{x})[X,T_1,T_2].
  \]
  Then the coefficients of $F$ (considered as a polynomial in $X$,
  $T_1$, $T_2$) form a generating set for the invariant field
  $K(\underline{x})_0^{\Sigma_n \times \PGL_3(K)}$.
\end{theorem}

\begin{proof}
  We may assume $n \ge 5$, since for $n \le 4$ all invariants of
  $\PGL_3(K)$ are constant, and the polynomial $F$ is the empty
  product, so we are claiming $K(\underline{x})_0^{\Sigma_n \times
  \PGL_3(K)} = K$ in this case, which is true.
  
  Write $L$ for the field extension of $K$ generated by the
  coefficients of $F$, and set
  \[ 
  {\mathcal M} := \left\{S \subseteq \{1 \upto n\} \mid |S| =
  5\right\}.
  \]
  Since the coefficients of $F$ are rational functions in the
  $c_{i,j,k,l,m}$, it follows that all elements from $L$ are
  $\PGL_3(K)$-invariant. Moreover, any $\pi \in \Sigma_n$ affords a
  permutation of $\mathcal M$, hence the product $F$, and therefore
  its coefficients, are fixed by~$\pi$. It follows that $L \subseteq
  K(\underline{x})_0^{\Sigma_n \times \PGL_3(K)}$.

  To prove the reverse inclusion, set
  \[
  {\mathcal C} := \left\{c_{i,j,k,l,m} \mid i,j,k,l,m \in \{1 \upto
  n\} \ \text{pairwise distinct}\right\},
  \]
  and for $S \in \mathcal M$ set
  \[
  {\mathcal C}_S := \left\{c_{i,j,k,l,m} \in {\mathcal C} \mid
  \{i,j,k,l,m\} = S\right\},
  \]
  so $\mathcal C$ is the disjoint union of all the ${\mathcal
  C}_S$. For $S \in \mathcal M$, the polynomial
  \[
  f_S := \prod_{c_{i,j,k,l,m} \in {\mathcal C}_S} (X - c_{i,j,k,l,m})
  \]
  has coefficients which are invariant under all permutations of the
  set $S$. With $S_0 := \{1,2,3,4,5\}$, \tref{2tN5} may be restated as
  \begin{equation} \label{3eqN5}
    K({\mathcal C}_{S_0})^{\Sigma_5} = K(x_{\nu,\mu} \mid \nu \in \{1
    \upto 5\},\mu = 0,1,2)^{\Sigma_5 \times \PGL_3(K)} =
    K(a_{S_0},b_{S_0})
  \end{equation}
  (where \tref{1tCs} was used for the first equality), so we obtain
  $f_{S_0} \in K(a_{S_0},b_{S_0})[X]$. Thus we can write $f_{S_0} =
  R(a_{S_0},b_{S_0},X)$, where $R$ is a rational function of three
  arguments (with the third argument not appearing in the denominator
  of $R$). But exactly the same will be true if we replace the indices
  $1,2,3,4,5$ by indices $i,j,k,l,m$ with $\{i,j,k,l,m\} = S$. So we
  obtain
  \begin{equation} \label{3eqFC}
    f_S = R(a_S,b_S,X)
  \end{equation}
  for all $S \in \mathcal M$ with $R \in K(Y,Z,X)$ a rational function
  not depending on $S$. This equation will be used later in the
  proof. Here we conclude that
  \[
  f := \prod_{c_{i,j,k,l,m} \in {\mathcal C}} (X - c_{i,j,k,l,m}) =
  \prod_{S \in {\mathcal M}} f_S \in K\left(a_S,b_S \mid S \in
  {\mathcal M}\right)[X].
  \]
  Let $\sigma$ be a $K$-automorphism of $K\left(a_S,b_S \mid S \in
  {\mathcal M}\right)$ which is given by a permutation $\psi$ of the
  set $\mathcal M$. Then by~\eqref{3eqFC}, $\sigma$ permutes the
  factors $f_S$ of~$f$ and therefore fixes~$f$. Thus the coefficients
  of~$f$ lie in the fixed field of all automorphisms~$\sigma$ of this
  type. Moreover, the $a_S$ and $b_S$ are pairwise distinct, since
  $a_S$ and $b_S$ are distinct, and for different sets $S$ they
  involve different sets of variables $x_{\nu,\mu}$. Hence we can use
  \lref{3lGalois}, which tells us that the coefficients of~$f$ lie in
  $L$. It follows that the field $K({\mathcal C})$ generated by the
  roots of~$f$ is the splitting field of~$f$ over $L$. Since the
  $c_{i,j,k,l,m} \in \mathcal C$ are pairwise distinct (as we defined
  $\mathcal C$ as a {\em set}), $f$ is separable, and therefore
  $K({\mathcal C})$ is Galois as a field extension of $L$. Assume that
  we can show that $\Gal\left(K({\mathcal C})/L\right)$ is contained
  in $\Sigma_n$ (i.e., every $\sigma$ in the Galois group is given by
  a permutation from $\Sigma_n$ acting on the $c_{i,j,k,l,m} \in
  \mathcal C$ by permuting the indices), then we have
  \[
  K(\underline{x})_0^{\Sigma_n \times \PGL_3(K)} = K({\mathcal
  C})^{\Sigma_n} \subseteq K({\mathcal C})^{\Gal\left(K({\mathcal
  C})/L\right)} = L
  \]
  (where \tref{1tCs}(a) was used for the first equation), and we are
  done. Thus all we need to show is
  \begin{equation} \label{3eqToShow}
    \Gal\left(K({\mathcal C})/L\right) \subseteq \Sigma_n.
  \end{equation}
  So take $\sigma \in \Gal\left(K({\mathcal C})/L\right)$. Since
  $K({\mathcal C})$ is the splitting field of~$f$ over $L$, $\sigma$
  permutes the set $\mathcal C$. Moreover, we have
  \[
  \prod_{S \in \mathcal M} \left(X - T_1 a_S - T_2 b_S\right) = F =
  \sigma(F) = \prod_{S \in \mathcal M} \left(X - T_1 \sigma(a_S) - T_2
  \sigma(b_S)\right).
  \]
  Since the roots of a polynomial are unique up to permutation, there
  exists a permutation $\psi$ of $\mathcal M$ such that
  \begin{equation} \label{3eqSigma}
    \sigma(a_S) = a_{\psi(S)} \quad \text{and} \quad \sigma(b_S) =
    b_{\psi(S)}
  \end{equation}
  for all $S \in \mathcal M$. Together with~\eqref{3eqFC}, this
  implies $\sigma(f_S) = f_{\psi(S)}$. Using the definition of $f_S$,
  this means that
  \[
  \prod_{c_{i,j,k,l,m} \in {\mathcal C}_S} \left(X -
  \sigma(c_{i,j,k,l,m})\right) = \prod_{c_{i,j,k,l,m} \in {\mathcal
  C}_{\psi(S)}} \left(X - c_{i,j,k,l,m}\right),
  \]
  so $\sigma({\mathcal C}_S) = {\mathcal C}_{\psi(S)}$.

  Fix an $S \in \mathcal M$ and pick a bijection
  $\map{\pi_0}{\psi(S)}{S}$. Define a $K$-automorphism
  \[
  \map{\phi_{\pi_0}}{K\left(x_{\nu,\mu} \mid \nu \in \psi(S),\mu =
  0,1,2\right)}{K\left(x_{\nu,\mu} \mid \nu \in S,\mu = 0,1,2\right)}
  \]
  by setting $\phi_{\pi_0}(x_{\nu,\mu}) := x_{\pi_0(\nu),\mu}$. Then
  for $\{i,j,k,l,m\} = \psi(S)$ we have $\phi_{\pi_0}(c_{i,j,k,l,m}) =
  \linebreak c_{\pi_0(i),\pi_0(j),\pi_0(k),\pi_0(l),\pi_0(m)}$, so
  $\phi_{\pi_0}(a_{\psi(S)}) = a_S$ and $\phi_{\pi_0}(b_{\psi(S)}) =
  b_S$. Together with~\eqref{3eqSigma} this implies
  $\left(\phi_{\pi_0} \circ \sigma\right)(a_S) = a_S$ and
  $\left(\phi_{\pi_0} \circ \sigma\right)(b_S) = b_S$. From
  $\sigma({\mathcal C}_S) = {\mathcal C}_{\psi(S)}$ we see that
  $\phi_{\pi_0} \circ \sigma$ maps $K({\mathcal C}_S)$ to itself.
  Therefore $\phi_{\pi_0} \circ \sigma$ restricted to $K({\mathcal
    C}_S)$ is a $K$-automorphism which fixes $a_S$ and $b_S$. But we
  have $K({\mathcal C}_S)^{\Sigma_S} = K(a_S,b_S)$, where $\Sigma_S$
  is the group of all permutations of $S$ (this is~\eqref{3eqN5}
  restated with the indices $1,2,3,4,5$ replaced by $i,j,k,l,m$ with
  $\{i,j,k,l,m\} = S$). By Galois theory, this implies that
  $\phi_{\pi_0} \circ \sigma$ restricted to $K({\mathcal C}_S)$ lies
  in $\Sigma_S$, i.e., there exists $\pi \in \Sigma_S$ such that
  $\left(\phi_{\pi_0} \circ \sigma\right)(c_{i,j,k,l,m}) =
  c_{\pi(i),\pi(j),\pi(k),\pi(l),\pi(m)}$ for all $i,j,k,l,m$ with
  $\{i,j,k,l,m\} = S$. Set $\map{\pi_S := \pi_0^{-1} \circ
    \pi}{S}{\psi(S)}$. Then $\pi_S$ is a bijection and we have
  $\sigma(c_{i,j,k,l,m}) =
  c_{\pi_S(i),\pi_S(j),\pi_S(k),\pi_S(l),\pi_S(m)}$ for all
  $i,j,k,l,m$ with $\{i,j,k,l,m\} = S$. This can be done with all $S
  \in \mathcal M$.

  In summary, we have a permutation $\psi$ of $\mathcal M$, and for
  each $S \in \mathcal M$ we have a bijection
  $\map{\pi_S}{S}{\phi(S)}$ such that
  \[
  \sigma(c_{i,j,k,l,m}) =
  c_{\pi_S(i),\pi_S(j),\pi_S(k),\pi_S(l),\pi_S(m)} \quad \text{where}
  \quad S = \{i,j,k,l,m\}.
  \]
  Being a field-automorphism, $\sigma$ preserves all algebraic
  relations that exist between the $c_{i,j,k,l,m}$. Thus we are
  exactly in the situation of \lref{3lPerm}, which tells us that
  $\sigma$ lies in $\Sigma_n$. Thus~\eqref{3eqToShow} is shown and the
  proof is complete.
\end{proof}

\begin{remark*}
  \begin{enumerate}
  \item Everything that was said in \rref{2rN5}(a) about the
    computation of the invariants~$a$ and~$b$ applies to the
    computation of the $a_S$ and $b_S$ used in \tref{3tNgen}, too. In
    particular, for each subset $S \subseteq \{1 \upto n\}$ with five
    elements, one only needs to evaluate two of the $c_{i,j,k,l,m}$ in
    order to calculate $a_S$ and $b_S$.
  \item As in the case of \tref{2tN5}, we can also deal with the case
    $\ch(K) = 2$ (see \rref{2rN5}(c)).
  \end{enumerate}
\end{remark*}

We will now turn to looking at separating properties of our
invariants. We need the following lemma.

\begin{lemma} \label{3lSeparating}
  Let $K$ be any field and let $g_1 \upto g_m \in K(x_1 \upto x_n)$ be
  rational functions in~$n$ indeterminates over $K$. Moreover, assume
  that $G$ is a finite group acting by $K$-automorphisms on the
  subfield $K(g_1 \upto g_m)$ generated by the $g_i$. Let $f_1 \upto
  f_r$ be generators of the invariant field, i.e., assume $K(g_1 \upto
  g_m)^G = K(f_1 \upto f_r)$. Then there exists a non-zero polynomial
  $h \in K[x_1 \upto x_n] \setminus \{0\}$ such that for all $\xi_1
  \upto \xi_n \in K$ with $h(\xi_1 \upto \xi_n) \ne 0$ the following
  holds: If $\eta_1 \upto \eta_n \in K$ are such that $f_i(\eta_1
  \upto \eta_n) = f_i(\xi_1 \upto \xi_n)$ for all $i \in \{1 \upto
  r\}$ (which is meant to imply that no zero-division occurs on either
  side of the equation), then there exists $\sigma \in G$ such that
  \[
  g_i(\eta_1 \upto \eta_n) = \left(\sigma(g_i)\right)(\xi_1 \upto
  \xi_n) \quad \text{for} \quad i \in \{1 \upto m\}.
  \]
  Moreover, $h$ can be chosen as the numerator of a polynomial in $f_1
  \upto f_n$ (viewed as a rational function in $K(x_1 \upto x_n)$).
\end{lemma}

\begin{proof}
  Parts of this proof are drawn from the proof of Theorem~3.9.13 in
  \mycite{Derksen:Kemper}. Take additional indeterminates $X$ and $T$,
  and form the polynomial
  \[
  F := \prod_{\sigma \in G} \left(X - \sum_{i=1}^m \sigma(g_i) \cdot
  T^{i-1}\right) \in K(x_1 \upto x_n)[X,T].
  \]
  $F$ is invariant under the action of $G$, thus all coefficients of
  $F$ lie in $K(g_1 \upto g_m)^G = K(f_1 \upto f_r)$. Let $c$ be a
  coefficient of $F$. Then we can write $c = F_c(f_1 \upto
  f_r)/H_c(f_1 \upto f_r)$ with $F_c,H_c \in K[T_1 \upto T_r]$
  polynomials and $H_c(f_1 \upto f_r) \ne 0$. Set $H \in K[T_1 \upto
  T_r]$ to be the lcm of all $H_c$ with~$c$ a coefficient of $F$. Thus
  $H(f_1 \upto f_r) \ne 0$. Let $h \in K[x_1 \upto x_n]$ be the
  numerator of $H(f_1 \upto f_r)$ (as a rational function in $K(x_1
  \upto x_n)$). Now assume we have $\xi_1 \upto \xi_n,\eta_1 \upto
  \eta_n \in K$ such that $h(\xi_1 \upto \xi_n) \ne 0$ and
  \begin{equation} \label{3eqFi}
      f_i(\xi_1 \upto \xi_n) = f_i(\eta_1 \upto \eta_n) \quad
      \text{for} \quad i \in \{1 \upto r\}.
  \end{equation}
  It follows that $\left(H(f_1 \upto f_r)\right)(\xi_1 \upto \xi_n)$
  is non-zero, and by~\eqref{3eqFi} the same is true for \linebreak
  $\left(H(f_1 \upto f_r)\right)(\eta_1 \upto \eta_n)$. Thus every
  coefficient~$c$ of $F$ can be evaluated at $(\xi_1 \upto \xi_n)$ and
  at $(\eta_1 \upto \eta_n)$, and we have $c(\xi_1 \upto \xi_n) =
  c(\eta_1 \upto \eta_n)$. For $\sigma \in G$, write $a_\sigma :=
  \sum_{i=1}^m \sigma(g_i) \cdot T^{i-1} \in K(x_1 \upto x_n)[T]$. It
  follows from the definition of $F$ that for every $\sigma \in G$ we
  have $F(a_\sigma) = 0$, where $X$ is taken as the main variable of
  $F$. Since $F$ is monic, it follows that an irreducible polynomial
  from $K[x_1 \upto x_n,T]$ which divides the denominator of
  $a_\sigma$ must also divide the denominator of at least one
  coefficient from $F$. Thus the fact that no zero-division occurs
  when substituting $(x_1 \upto x_n) = (\xi_1 \upto \xi_n)$ or $(x_1
  \upto x_n) = (\eta_1 \upto \eta_n)$ into the coefficients of $F$
  implies that also all $a_\sigma$ and hence all $\sigma(g_i)$ can be
  evaluated at $(\xi_1 \upto \xi_n)$ and at $(\eta_1 \upto \eta_n)$.
  Using $c(\xi_1 \upto \xi_n) = c(\eta_1 \upto \eta_n)$ for all
  coefficients~$c$ of $F$, we conclude that
  \[
  \prod_{\sigma \in G} \left(X - \sum_{i=1}^m
  \left(\sigma(g_i)\right)(\xi_1 \upto \xi_n) \cdot T^{i-1}\right) =
  \prod_{\sigma \in G} \left(X - \sum_{i=1}^m
  \left(\sigma(g_i)\right)(\eta_1 \upto \eta_n) \cdot T^{i-1}\right).
  \]
  The right hand side, regarded as a polynomial in $X$, has the zero
  $\sum_{i=1}^m g_i(\eta_1 \upto \eta_n) \cdot T^{i-1}$. This must
  also be a zero of the left hand side, hence there exists a $\sigma
  \in G$ such that
  \[
  \sum_{i=1}^m g_i(\eta_1 \upto \eta_n) \cdot T^{i-1} = \sum_{i=1}^m
  \left(\sigma(g_i)\right)(\xi_1 \upto \xi_n) \cdot T^{i-1}.
  \]
  Comparing coefficients in $T$ now yields $g_i(\eta_1 \upto \eta_n) =
  \left(\sigma(g_i)\right)(\xi_1 \upto \xi_n)$ for $i \in \{1 \upto
  m\}$, as desired.
\end{proof}

If we have~$n$ points $P_1 \upto P_n \in \PP^2(K)$ in projective
2-space such that no three of the $P_i$ are collinear, we can evaluate
the invariants $a_S$ and $b_S$ at $(P_1 \upto P_n)$ for every subset
$S \subseteq \{1 \upto n\}$ with $|S| = 5$. Thus for each $S$ we
obtain a vector $\left(a_S(\underline{P}),b_S(\underline{P})\right)
\in K^2$. We will consider the {\em distribution} of these vectors for
all subsets $S$. This distribution is adequately represented by the
polynomial
\[
F_{P_1 \upto P_n} := \prod
\begin{Sb}
  S \subseteq \{1 \upto n\}, \\
  |S| = 5
\end{Sb}
\left(X - T_1 a_S(P_1 \upto P_n) - T_2 b_S(P_1 \upto P_n)\strut\right)
\in K[X,T_1,T_2]
\]
with $X$, $T_1$, $T_2$ indeterminates. It is our goal to use these
distributions for two point configurations $(P_1 \upto P_n)$ and $(Q_1
\upto Q_n) \in \left(\PP^2(K)\right)^n$ to determine if $(P_1 \upto
P_n)$ can be transformed into $(Q_1 \upto Q_n)$ by a projective
transformation and a relabeling the points. We call a point
configuration $(P_1 \upto P_n) \in \left(\PP^2(K)\right)^n$
\df{reconstructible from the joint distribution of $a$'s and $b$'s} if
for any other $(Q_1 \upto Q_n) \in \left(\PP^2(K)\right)^n$ with
\[
F_{P_1 \upto P_n} = F_{Q_1 \upto Q_n}
\]
there exist a permutation $\pi \in \Sigma_n$ and a transformation $g
\in \PGL_3(K)$ such that
\[
Q_i = g\left(P_{\pi(i)}\right)
\]
for all $i \in \{1 \upto n\}$. In order to be able to apply
\tref{3tNgen}, we assume that the characteristic of $K$ is not~2.

\begin{cor} \label{3cSeparating}
  With the above notation there exists a non-zero polynomial $f \in
  K[\underline{x}] = K\left[x_{i,j} \mid \right. \linebreak \left. i =
    1 \upto n, j = 0,1,2\right]$ which for each~$i$ is homogeneous as
  a polynomial in $x_{i,0}$, $x_{i,1}$, $x_{i,2}$, such that every
  point configuration $(P_1 \upto P_n) \in \left(\PP^2(K)\right)^n$
  with $f(P_1 \upto P_n) \ne 0$ is reconstructible from the joint
  distribution of $a$'s and $b$'s.
\end{cor}

\begin{proof}
  Let $G := \Sigma_n$ be the symmetric group acting on the set
  \[
  {\mathcal C} := \left\{c_{i,j,k,l,m} | i,j,k,l,m \in \{1 \upto n\} \
  \text{pairwise distinct}\right\}
  \]
  by permuting the indices of the $c$'s. Thus $G$ acts by
  $K$-automorphisms on the field $K({\mathcal C})$ generated by the
  $c_{i,j,k,l,m}$. By \tref{1tCs}(a) we have that $K({\mathcal C}) =
  K(\underline{x})_0^{\PGL_3(K)}$. Write $f_1 \upto f_r \in
  K({\mathcal C})$ for the coefficients of the polynomial $F$ defined
  in \tref{3tNgen}. Then \tref{3tNgen} says that
  \[
  K({\mathcal C})^G = K(\underline{x})_0^{\Sigma_5 \times \PGL_3(K)} =
  K(f_1 \upto f_r).
  \]
  Thus we are exactly in the situation of \lref{3lSeparating}, which
  gives us a polynomial $h \in K[\underline{x}]$ with the properties
  stated in the lemma. Since $h$ is the numerator of a polynomial
  involving the $f_i$ (and therefore the $c_{i,j,k,l,m}$, which lie in
  $K(\underline{x})_0$), $f$ is homogeneous as a polynomial in
  $x_{i,0}$, $x_{i,1}$, $x_{i,2}$ for each~$i$ (see the proof of
  \tref{1tCs}). Let $f$ be the product of~$h$ and all determinants
  $[i,j,k]$ (defined before~\eqref{1eqC}) with $1 \le i < j < k \le
  n$. Now take $(P_1 \upto P_n) \in \left(\PP^2(K)\right)^n$ and
  assume $f(P_1 \upto P_n) \ne 0$. Moreover, take $(Q_1 \upto Q_n) \in
  \left(\PP^2(K)\right)^n$ with $F_{P_1 \upto P_n} = F_{Q_1 \upto
  Q_n}$. This means that all coefficients of $F$ take the same value
  when evaluated at $(P_1 \upto P_n)$ or at $(Q_1 \upto Q_n)$, so
  $f_i(P_1 \upto P_n) = f_i(Q_1 \upto Q_n)$ for $i = 1 \upto r$. By
  \lref{3lSeparating} there exists a $\pi \in G$ such that
  \[
  c_{i,j,k,l,m}(Q_1 \upto Q_n) = \left(\pi(c_{i,j,k,l,m})\right)(Q_1
  \upto P_n) = c_{i,j,k,l,m}(P_{\pi(1)} \upto P_{\pi(n)})
  \]
  for all $i,j,k,l,m \in \{1 \upto n\}$ pairwise distinct. Since
  $f(P_1 \upto P_n) \ne 0$ guarantees that no three of the $P_i$ are
  collinear, it follows from \tref{1tCs}(b) that there exists a $g \in
  \PGL_3(K)$ such that $Q_i = g\left(P_{\pi(i)}\right)$ for all $i = 1
  \upto n$. So $(P_1 \upto P_n)$ is reconstructible from the joint
  distribution of $a$'s and $b$'s.
\end{proof}

\section{Other groups} \label{4sGroups}

In this paper and in~[\citenumber{Boutin.Kemper}], we only considered
some very specific (though important) groups, namely projective,
Euclidean and volume-preserving groups. In this section we will look
at more general groups. The goal is to use reconstruction theorems
such as \cref{3cSeparating} for deriving reconstructibility statements
which classify a point configuration modulo any subgroup of the
original group. We will be more precise after proving the following
lemma.

\begin{lemma} \label{4lSn}
  Let~$n$ and~$m$ be integers with $0 < m < n$. Then the natural
  action of the symmetric group $\Sigma_n$ on the set
  \[
  {\mathfrak X} := \left\{M \subseteq \{1 \upto n\} \mid |M| = m\right\}
  \]
  is faithful.
\end{lemma}

\begin{proof}
  Suppose that for a $\pi \in \Sigma_n$ we have $\pi(M) = M$ for all
  $M \in {\mathfrak X}$. Take any $i \in \{1 \upto n\}$. Then
  \[
  \pi\left(\{i\}\right) = \pi\left(\strut\right.\bigcap
  \begin{Sb}
    M \in {\mathfrak X}, \\
    i \in M
  \end{Sb}
  M\left.\strut\right) = \bigcap
  \begin{Sb}
    M \in {\mathfrak X}, \\
    i \in M
  \end{Sb}
  \pi(M) = \bigcap
  \begin{Sb}
    M \in {\mathfrak X}, \\
    i \in M
  \end{Sb}
  M = \{i\},
  \]
  where the second equality follows from the injectiveness of~$\pi$.
  Hence $\pi(i) = i$.
\end{proof}

Let $X$ be any set (e.g., a projective or linear space) and let $G$ be
a group acting on $X$. For $M \subseteq X$ and $g \in G$ we write
\[
g(M) := \{g(x) \mid x \in M\} \quad \text{and} \quad G(M) := \{g(M)
\mid g \in G\}.
\]
Thus $G(M)$ is a subset of the power set ${\mathfrak P}(X)$ of $X$. We
may think of $X$ as a set of points and of $M$ (if finite) as a point
configuration, where the labeling of the points in $M$ is already
disregarded since we are considering $M$ as a set. Then $G(M)$ is the
class of all point-configurations which are ``congruent'' to $M$,
where the concept of ``congruence'' is given by the $G$-action. Fix a
positive integer~$m$. For $C \subseteq X$ a finite subset let
$\mu_{m,G}(C)$ be the multiset formed of all $G(M)$ with $M \subseteq
C$ and $|M| = m$. Formally, $\mu_{m,G}(C)$ may be defined as the
function ${\mathfrak P}\left({\mathfrak P}(X)\right) \to \ZZ$
assigning to each subset ${\mathfrak X} \subseteq {\mathfrak P}(X)$
the number $\left|\left\{M \subseteq C \mid |M| = m,\ G(M) =
    {\mathfrak X}\right\}\right|$. So $\mu_{m,G}(C)$ may be viewed as
the distribution of all $m$-subsets of $C$ up to the $G$-action.
Clearly for any $g \in G$ we have $\mu_{m,G}\left(g(C)\right) =
\mu_{m,G}(C)$. We call $C$ \df{reconstructible from $m$-subsets modulo
  $G$} if for every finite subset $D \subseteq X$ with $\mu_{m,G}(D) =
\mu_{m,G}(C)$ there exists $g \in G$ with $D = g(C)$.

In this language, \cref{3cSeparating} implies that ``almost'' all
finite subsets of $\PP^2(K)$ are reconstructible from 5-subsets modulo
$\PGL_3(K)$. Likewise, Theorem~1.6 from \mycite{Boutin.Kemper} says
that almost all finite subsets of $K^m$ (of size $\ge m + 2$) are
reconstructible from 2-subsets modulo the Euclidean group $\AO_m$.

\begin{theorem} \label{4tSubgroup}
  With the above notation assume that
  \begin{enumerate}
  \renewcommand{\theenumi}{\roman{enumi}}
  \item \label{4i1} $C$ is reconstructible from $m$-subsets modulo $G$,
  \item \label{4i2} for $M,N \subseteq C$ with $|M| = |N| = m$, we
    have that $G(M) = G(N)$ implies $M = N$, and
  \item \label{4i3} there exists a subset $\hat{M} \subseteq C$ with
    $|\hat{M}| = m + 1$ such that
    \[
    \left\{g \in G | g(x) = x \ \text{for all} \ x \in \hat{M}\right\}
    = \{\id\}.
    \]
  \end{enumerate}
  Then for every subgroup $H \le G$, $C$ is reconstructible from $(m +
  1)$-subsets modulo $H$.
\end{theorem}

\begin{proof}
  Let $D \subseteq X$ be a finite set with
  \begin{equation} \label{4eqHyp}
    \mu_{m+1,H}(D) = \mu_{m+1,H}(C).
  \end{equation}
  We wish to show that there exists $h \in H$ with $D = h(C)$. Since
  $\left|\mu_{m+1,H}(C)\right| = \binom{|C|}{m+1}$ and $|C| \ge m+1$
  by the assumption~\eqref{4i3}, \eqref{4eqHyp} certainly implies $|D|
  = |C|$. Take any subset $M \subseteq C$ with $|M| = m$. The
  assumption~\eqref{4i3} implies that $|C| > m$, so there exists $x
  \in C \setminus M$. Set $M' := M \cup \{x\}$. By~\eqref{4eqHyp}
  there exists $N' \subseteq D$ with $|N'| = m+1$ such that $H(M') =
  H(N')$. So there exists $g \in H$ with $M' = g(N')$. Thus we have $M
  \subset g(N')$, so there exists a subset $N \subset N'$ with $|N| =
  m$ and $M = g(N)$. This implies $G(M) = G(N)$. Since $M$ was taken
  to be an arbitrary $m$-subset of $C$, it follows that $\mu_{m,G}(C)
  \subseteq \mu_{m,G}(D)$ (observe that by~\eqref{4i2} the multiset
  $\mu_{m,G}(C)$ has no multiplicities). Since $|C| = |D|$, the
  cardinalities of $\mu_{m,G}(C)$ and $\mu_{m,G}(D)$ also coincide,
  and we conclude $\mu_{m,G}(C) = \mu_{m,G}(D)$. Note that this
  implies that the assumption~\eqref{4i2} also holds for $C$ replaced
  by $D$. But the main consequence of $\mu_{m,G}(C) = \mu_{m,G}(D)$ is
  that by~\eqref{4i1} there exists $g \in G$ such that
  \begin{equation} \label{4eqG}
    g(C) = D.
  \end{equation}
  Now we consider the subset $\hat{M} \subseteq C$ given
  by~\eqref{4i3}. By~\eqref{4eqHyp} we have a subset $\hat{N}
  \subseteq D$ with $|\hat{N}| = m + 1$ and $H(\hat{M}) = H(\hat{N})$.
  So there exists $h \in H$ with $h(\hat{M}) = \hat{N}$. Take any $M
  \subset \hat{M}$ with $|M| = m$. Then
  \[
  N := h(M) \subset h(\hat{M}) = \hat{N} \subseteq D.
  \]
  $N = h(M)$ implies $G(N) = G(M)$. For $\tilde{N} := g(M)$ (with~$g$
  from \eqref{4eqG}) we also have $G(\tilde{N}) = G(M)$, so
  $G(\tilde{N}) = G(N)$. By~\eqref{4eqG}, $\tilde{N} \subseteq D$, and
  since~\eqref{4i2} also holds with $C$ replaced by $D$, we conclude
  that $\tilde{N} = N$, e.i., $g(M) = h(M)$. This holds for any
  $m$-subset $M \subset \hat{M}$. Thus
  \[
  g(\hat{M}) = g\left(\strut\right.\bigcup
  \begin{Sb}
    M \subset \hat{M}, \\
    |M| = m
  \end{Sb}
  M\left.\strut\right) = \bigcup
  \begin{Sb}
    M \subset \hat{M}, \\
    |M| = m
  \end{Sb}
  g(M) = \bigcup
  \begin{Sb}
    M \subset \hat{M}, \\
    |M| = m
  \end{Sb}
  h(M) = h\left(\strut\right.\bigcup
  \begin{Sb}
    M \subset \hat{M}, \\
    |M| = m
  \end{Sb}
  M\left.\strut\right) = h(\hat{M}).
  \]
  It follows that $h^{-1} \circ g$ restricts to a permutation~$\pi$ of
  $\hat{M}$. Since $g(M) = h(M)$ for all $m$-subsets $M \subset
  \hat{M}$, $\pi(M) = M$ for all these $M$. It follows by \lref{4lSn}
  that $\pi = \id$. Thus $g|_{_{\hat{M}}} = h|_{_{\hat{M}}}$ (the
  restrictions to $\hat{M}$ coincide). Now~\eqref{4i3} yields $g = h$,
  so~\eqref{4eqG} implies $h(C) = D$, which completes the proof.
\end{proof}

\begin{cor}[consequence of \tref{4tSubgroup} and \cref{3cSeparating}]
  \label{4cSubPGL}
  Let $K$ be an infinite field and $n \ge 6$ an integer. Then there
  exists a non-zero polynomial $f \in K[\underline{x}] =
  K\left[x_{i,j} | \mid i = 1 \upto n, j = 0,1,2\right]$ which for
  each~$i$ is homogeneous as a polynomial in $x_{i,0}$, $x_{i,1}$,
  $x_{i,2}$, such that for every point configuration $(P_1 \upto P_n)
  \in \left(\PP^2(K)\right)^n$ with $f(P_1 \upto P_n) \ne 0$ the set
  $\{P_1 \upto P_n\}$ is reconstructible from 6-subsets modulo $G$ for
  every subgroup $G \le \PGL_3(K)$.
\end{cor}

\begin{proof}
  By \cref{3cSeparating} there exists a non-zero polynomial
  $\tilde{f}$ such that all $(P_1 \upto P_n)$ with $\tilde{f}(P_1
  \upto P_n) \ne 0$ is reconstructible from the joint distribution of
  $a$'s and $b$'s. In particular, this means that for such $(P_1 \upto
  P_n)$ the set $\{P_1 \upto P_n\}$ is reconstructible from 5-subsets
  modulo $\PGL_3$. This provides the hypothesis~\eqref{4i1} of
  \tref{4tSubgroup}. The hypothesis~\eqref{4i2} can also be turned
  into an open condition on $(P_1 \upto P_n)$. Indeed, it is enough to
  impose that for distinct 5-subsets $M$ and $N$ of $\{P_1 \upto
  P_n\}$, the pairs $\left(a(M),b(M)\right)$ and
  $\left(a(N),b(N)\right)$ (with $a$ and $b$ the ($\Sigma_5 \times
  \PGL_3$)-invariants defined in \tref{2tN5}) are also distinct. To
  make sure that~\eqref{4i3} also holds, it suffices by the uniqueness
  statement in \lref{1lPGLOrbits} that there exist four points in
  $\{P_1 \upto P_n\}$ such that no three of them are collinear, which
  is also an open condition. Finally, one should impose the condition
  that the $P_i$ are pairwise distinct to ensure that the set of the
  $P_i$ really has size~$n$.
\end{proof}

In the following corollary, $K$ is any field and $V$ is an
$m$-dimensional vector space over $K$. We write $V^n$ for the direct
sum of~$n$ copies of $V$, and $K[V^n]$ for the ring of polynomials on
$V^n$. $\ASL^\pm(V)$ is the group generated by all linear
transformations of $V$ with determinant $\pm 1$ and all translations
of $V$.

\begin{cor}[consequence of \tref{4tSubgroup} and Theorem~3.7 from
  \mycite{Boutin.Kemper}] \label{4cSubAGL} \mbox{} \linebreak Assume
  $n \ge m + 2$.  Then there exists a non-zero polynomial $f \in
  K[V^n]$ such that for $(P_1 \upto P_n) \in V^n$ with $f(P_1 \upto
  P_n) \ne 0$, the set $\{P_1 \upto P_n \}$ is reconstructible from
  $(m+2)$-subsets modulo $G$ for every subgroup $G \le \ASL^\pm(V)$.
\end{cor}

\begin{proof}
  Theorem~3.7 from \mycite{Boutin.Kemper} says that there exists
  $\tilde{f} \in K[V^n] \setminus \{0\}$ such that all $(P_1 \upto
  P_n) \in V^n$ with $\tilde{f}(P_1 \upto f_n) \ne 0$ are
  reconstructible (up to the actions of $\ASL^\pm(V)$ and the
  symmetric group $\Sigma_n$) from the distribution of volumes of
  parallelepiped spanned by $(m+1)$-subsets. In particular, for these
  $(P_1 \upto P_n)$, the set $\{P_1 \upto P_n\}$ is reconstructible
  from $(m+1)$-subsets modulo $\ASL^\pm(V)$. Moreover, imposing that
  for distinct $(m+1)$-subsets of $\{P_1 \upto P_n\}$ the volumes of
  the parallelepiped spanned by these subsets also differ is an open
  condition. Finally, the assumption~\eqref{4i3} in \tref{4tSubgroup}
  is satisfied if $\{P_1 \upto P_n\}$ contains $m+1$ points which span
  a parallelepiped of non-zero volume.
\end{proof}

\begin{remark} \label{4rInvariants}
  Suppose that in the situation of \cref{4cSubAGL} we have (rational)
  invariants $f_1 \upto f_r \linebreak \in K(V^{m+2})^{\Sigma_{m+2}
    \times G}$ (where $K(V^{m+2})$ is the rational function field on
  $V^{m+2}$ and $G$ is the subgroup of $\ASL^\pm(V)$ which is
  considered) such that for a non-empty Zariski-open subset $S
  \subseteq V^{m+2}$ the invariants $f_i$ can be evaluated on $S$, and
  for $(P_1 \upto P_{m+2}) \in S$ and $(Q_1 \upto Q_{m+2}) \in
  V^{m+2}$ we have that $f_i(P_1 \upto P_{m+2}) = f_i(Q_1 \upto
  Q_{m+2})$ for all~$i$ implies that $Q_i = g(P_{\pi(i)})$ with $g \in
  G$ and $\pi \in \Sigma_{m+2}$. Then it follows from \cref{4cSubAGL}
  that for $n \ge m + 2$ there exists $f \in K[V^n] \setminus \{0\}$
  such that all $(P_1 \upto P_n) \in V^n$ with $f(P_1 \upto P_n) \ne
  0$ are reconstructible (modulo the actions of $G$ and $\Sigma_n$)
  from the joint distribution of $f_1 \upto f_r$ (i.e., the
  distribution of the values $(f_1(M) \upto f_r(M)) \in K^r$, where
  $M$ ranges through all $(m + 2)$-subsets of $\{P_1 \upto P_n\}$).

  The analogous remark applies in the situation of \cref{4cSubPGL}.
\end{remark}

\begin{ex} \label{4exCounter}
  This example shows that in \cref{4cSubAGL} the number $m + 2$ cannot
  be reduced to a lower number. Consider the case $m = 1$ (i.e., $V =
  K$, and let $G \cong K$ be the group of all translations. Consider a
  point-configuration $C := \{P_1 \upto P_n\} \subseteq K$ and its
  negative $-C := \{-P_1 \upto -P_n\}$. For a 2-subset $\{P_i,P_j\}
  \subseteq C$ the group element $g := -P_i - P_j$ yields
  \[
  g(\{P_i,P_j\}) = \{-P_j,-P_i\} \subseteq -C.
  \]
  Hence $\mu_{2,G}(C) = \mu_{2,G}(-C)$. But clearly $C$ and $-C$ are
  only congruent modulo $G$ if $C$ has a special symmetry property.
  Thus there exists no non-empty Zariski-open subset $S \subseteq K^n$
  such that all $n$-subsets of $K$ formed from tuples from $S$ are
  reconstructible from 2-subsets modulo $G$. This example shows that
  also in \tref{4tSubgroup} the number $m + 1$ cannot be decreased.
  
  However, by \cref{4cSubAGL}, for every $n \ge 3$ there exists an $f
  \in K[x_1 \upto x_n] \setminus \{0\}$ such that for $(P_1 \upto P_n)
  \in K^n$ with $f(P_1 \upto P_n) \ne 0$, the set $\{P_1 \upto P_n\}$
  is reconstructible from 3-subsets modulo $G$. We can also give
  invariants in $K[x_1,x_2,x_3]^{\Sigma_3 \times G}$ as in
  \rref{4rInvariants}, which can be found easily by using the
  invariant theory package in Magma~[\citenumber{magma}]. They are
  \[
  f_1 = x^2 + y^2 + z^2 - x y - x z - y z, \quad f_2 = (2 x - y - z)
  (2 y - x - z) (2 z - x - y),
  \]
  so almost all $n$-point configurations are determined up to
  $\Sigma_n \times G$ by the distribution of the vectors
  $\left(f_1(x_i,x_j,x_k),f_2(x_i,x_j,x_k)\right)$ for $\{i,j,k\}
  \subseteq \{1 \upto n\}$.
\end{ex}

\bibliographystyle{mybibstyle} \bibliography{bib}

\bigskip

\begin{center}
\begin{tabular}{lll}
  Mireille Boutin & & Gregor Kemper \\
  Department of Mathematics  & & Technische Universit\"at M\"unchen \\
  Purdue University & & Zentrum Mathematik - M11 \\
  150 N. University St.  & & Boltzmannstr. 3 \\
  West Lafayette, IN 47907 & & 85\,748 Garching \\
  USA & & Germany \\
  {\tt boutin$@$math.purdue.edu} & & {\tt kemper$@$ma.tum.de}
\end{tabular}
\end{center}

\end{document}